\documentclass[aop,MSNbibl,nameyear,dvips]{arximspdf}

\doi{10.1214/11-AOP694}
\volume{41}
\issue{2}
\pubyear{2013}
\firstpage{914}
\lastpage{960}

\makeatletter

\newtheorem{theorem}{Theorem}
\newtheorem{corollary}[theorem]{Corollary}
\newproclaim{definition}[theorem]{Definition}
\newtheorem{lemma}[theorem]{Lemma}
\newproclaim{notation}[theorem]{Notation}
\newtheorem{proposition}[theorem]{Proposition}
\newproclaim{remark}[theorem]{Remark}
\newproclaim{comment}[theorem]{Comment}
\newproclaim{Assumption}{Assumption}

\newcommand{\iint}{\int\!\!\int}

\makeatother

\begin{document}
\begin{frontmatter}

\title{Rosenthal-type inequalities for the maximum of partial sums of
stationary processes and examples}
\runtitle{Rosenthal-type inequalities for stationary processes}

\begin{aug}
\author[A]{\fnms{Florence} \snm{Merlev\`{e}de}\ead[label=e1]{florence.merlevede@univ-mlv.fr}}
\and
\author[B]{\fnms{Magda} \snm{Peligrad}\corref{}\thanksref{t1}\ead[label=e2]{peligrm@ucmail.uc.edu}}
\runauthor{F. Merlev\`{e}de and M. Peligrad}
\affiliation{Universit\'{e} Paris-Est and University of Cincinnati}
\address[A]{LAMA and UMR 8050 CNRS\\
Universit\'{e} Paris Est-Marne la Vall\'{e}e\\
5 Boulevard Descartes 77 454 Marne La Vall\'{e}e\\
France\\
\printead{e1}} 
\address[B]{Department of Mathematical Sciences\\
University of Cincinnati\\
PO Box 210025\\
Cincinnati, Ohio 45221-0025\\
USA\\
\printead{e2}}
\end{aug}

\thankstext{t1}{Supported in part by a Charles Phelps Taft Memorial
Fund grant, and
NSA Grants H98230-09-1-0005 and H98230-11-1-0135.}

\received{\smonth{9} \syear{2010}}
\revised{\smonth{5} \syear{2011}}

%
\begin{abstract}
The aim of this paper is to propose new Rosenthal-type inequalities for
moments of order higher than $2$ of the maximum of partial sums of
stationary sequences including martingales and their generalizations.
As in the recent results by Peligrad et al. [\textit{Proc. Amer. Math.
Soc.} \textbf{135} (2007) 541--550] and Rio [\textit{J. Theoret.
Probab.} \textbf{22} (2009) 146--163], the estimates of the moments are
expressed in terms of the norms of projections of partial sums. The
proofs of the results are essentially based on a new maximal inequality
generalizing the Doob maximal inequality for martingales and dyadic
induction. Various applications are also provided.
\end{abstract}

%

\end{frontmatter}

\section{Introduction}

For independent random variables, the Rosenthal inequalities relate
moments of
order higher than $2$ of partial sums of random variables to the
variance of
partial sums. One variant of this inequality is the following [see
\citet{Ros70}, page 279]: let $(X_{k})_{k}$ be independent and
centered real valued
random variables with finite moments of order $p$, $p\geq2$. Then for every
positive integer $n$,
%
%
\begin{equation}\label{Rosenind}%
{\mathbf{E}}\Bigl({\max_{1\leq j\leq n}} |S_{j}|^{p}\Bigr)\ll\sum_{k=1}%
^{n}{\mathbf{E}}(|X_{k}|^{p})+ \Biggl(\sum_{k=1}^{n}{\mathbf{E}}(X_{k}%
^{2}) \Biggr)^{p/2},
\end{equation}
where $S_{j}=\sum_{k=1}^{j}X_{k}$. Unless otherwise specified,
throughout the
paper the notation $a_{n}\ll b_{n}$ means that there exists a numerical
constant $C_{p}$ depending only on $p$ (and not on the underlying random
variables and neither on $n$) such that $a_{n}\leq C_{p}b_{n}$, for all
positive integers $n$.

Besides being useful to compare the norms ${\mathbf{L}}^{p}$ and
${\mathbf{L}}^{2}$ of partial sums, these inequalities are important
tools for
obtaining a variety of results, including tightness of the empirical process
[see the proof of Theorem 22.1 in \citet{Bil68}], convergence
rates with
respect to the strong law of large numbers [see, e.g.,
\citet{Wit85}]
or almost sure invariance principles; see \citet{Wu07} and
\citet{Gou10} for
recent results. Since the 1970s, there has been a great amount of work
which has
extended inequality (\ref{Rosenind}) to dependent sequences. See, for
instance, among many others: \citet{Pel85} and \citet{Sha95}
for the case of
$\rho$-mixing sequences; \citet{Sha88}, \citet{Pel89} and
\citet{Ute91} for the
case of $\phi$-mixing sequences; \citet{PelGut99} and \citet
{UtePel03} for interlaced mixing; Theorem 2.2 in \citet{Vie97} for
$\beta$-mixing
processes; Theorem 6.3 in \citet{Rio00} for the strongly mixing
case; \citet{Ded01} and \citet{Rio09} for projective criteria.

The main goal of this paper is to generalize the Rosenthal inequality from
sequences of independent variables to stationary dependent sequences including
martingales, allowing then to consider examples that are not necessarily
dependent in the sense of the dependence structures mentioned above.

In order to present our results, let us first introduce some notations and
definitions used throughout the paper.
\begin{notation}
\label{notation1} Let $(\Omega,{\mathcal{A}},{\mathbf{P}})$ be a probability
space, and let $T\dvtx\Omega\mapsto\Omega$ be a bijective bi-measurable
transformation preserving the probability ${\mathbf{P}}$. Let $\mathcal
{F}%
_{0}$ be a $\sigma$-algebra of $\mathcal{A}$ satisfying $\mathcal{F}%
_{0}\subseteq T^{-1}(\mathcal{F}_{0})$. We then define the nondecreasing
filtration $(\mathcal{F}_{i})_{i\in\mathbb{Z}}$ by $\mathcal{F}_{i}%
=T^{-i}(\mathcal{F}_{0})$ and the stationary sequence $(X_{i})_{i\in
\mathbb{Z}}$ by $X_{i}=X_{0}\circ T^{i}$, where $X_{0}$ is a real-valued
random variable. The sequence will be called adapted to the filtration
$(\mathcal{F}_{i})_{i\in\mathbb{Z}}$ if $X_{0}$ is ${}\mathcal{F}_{0}%
$-measurable. The following notations will also be used: ${\mathbf{E}}%
_{k}(X)={\mathbf{E}}(X|\mathcal{F}_{k})$ and the norm in ${\mathbf
{L}}^{p}$ of
$X$ is denoted by $\|X\|_{p}$. Let $S_{n}=\sum_{j=1}^{n}X_{j}$.
\end{notation}

In the rest of this section the sequence $(X_{i})_{i\in\mathbb{Z}}$ is
assumed to be stationary and adapted to $(\mathcal{F}_{i})_{i\in\mathbb{Z}}$
and the variables are in ${\mathbf{L}}^{p}$.

If $(X_{k})_{k}$ are stationary martingale differences, the martingale
form of inequality (\ref{Rosenind}) is
%
%
\begin{equation}\label{Rosenmart}\quad
\Bigl\|{\max_{1\leq j\leq n}} |S_{j}|\Bigr\|_{p}\ll n^{1/p}%
\|X_{1}\|_{p}+ \Biggl\Vert\sum_{k=1}^{n}{\mathbf{E}}_{k-1}(X_{k}%
^{2}) \Biggr\Vert_{p/2}^{1/2} \qquad\mbox{for any }p\geq2;
\end{equation}
see \citet{Bur73}. One of our goals is to replace the last term in this
inequality with a new one containing terms of the form $\Vert{\mathbf
{E}}%
_{0}(S_{n}^{2})\Vert_{p/2}$. The reason for introducing this term comes from
the fact that for many stationary sequences $\Vert{\mathbf{E}}_{0}(S_{n}
^{2})\Vert_{p/2}$ is closer to the variance of partial sums. In
addition, we
are interested in pointing out a Rosenthal-type inequality for a larger
class of
stationary adapted sequences that includes the martingale differences
as a
special case.

Two recent results by \citet{PelUte05} and \citet{WuZha08} show
that
\[
\Bigl\|{\max_{1\leq j\leq n}} |S_{j}|\Bigr\|_{p}\ll n^{1/p}\Biggl(
\|X_{1}\|_{p}+\sum_{k=1}^{n}\frac{1}{k^{1+1/p}}\|{\mathbf{E}}%
_{0}(S_{k})\|_{p}\Biggr) \qquad\mbox{for any }1\leq p\leq2.
\]

To find a suitable extension of this inequality for $p>2$, the first
step in
our approach is to establish the following maximal inequality that has
interest in itself:%
%
%
\begin{eqnarray}\label{max}%
\Bigl\Vert{\max_{1\leq j\leq n}}|S_{j}|\Bigr\Vert_{p}\ll n^{1/p}\Biggl(
{\max_{1\leq j\leq n}}\Vert S_{j}\Vert_{p}/n^{1/p}+\sum_{k=1}^{n}\frac
{1}{k^{1+1/p}}\|{\mathbf{E}}_{0}(S_{k})\|_{p}\Biggr)
\nonumber\\[-8pt]\\[-8pt]
&&\eqntext{\mbox{for any }p>1.}
\end{eqnarray}
This inequality can be viewed as generalization of the well-known Doob
maximal inequality for martingales. For a more precise version than
(\ref{max}), with constants specified, see our inequality (\ref{cons1sub}).

Then we combine inequality (\ref{max}) with several inequalities for
$\Vert S_{n}\Vert_{p}$ that will further be established in this paper.

As we shall see in Section~\ref{sectiongeneralpowers}, by a direct approach
using dyadic induction combined with maximal inequality (\ref{max}), we
shall prove that, for any $p>2$,
%
%
\begin{eqnarray}\label{inedirect}%
\Bigl\Vert{\max_{1\leq j\leq n}} |S_{j}|\Bigr\Vert_{p}
&\ll& n^{1/p}\Biggl(
\|X_{1}\|_{p}+\sum_{k=1}^{n}\frac{1}{k^{1+1/p}}\|{\mathbf{E}}%
_{0}(S_{k})\|_{p}\nonumber\\[-8pt]\\[-8pt]
&&\hphantom{n^{1/p}\Biggl(}
{}+ \Biggl(\sum_{k=1}^{n}\frac{1}{k^{1+2\delta/p}}\Vert
{\mathbf{E}}_{0}(S_{k}^{2})\Vert_{p/2}^{\delta} \Biggr)^{1/(2\delta)}\Biggr)
,\nonumber
\end{eqnarray}
where $\delta=\min(1,1/(p-2))$. For $2<p\leq3$ our inequality provides a
maximal form for Theorem 3.1 in \citet{Rio09}. When $p\geq4$, we
shall see that
the last term in the right-hand side dominates the second term, so that the
second term can be omitted in this case. Inequality (\ref{inedirect}) shows
that in order to relate $\Vert{\max_{1\leq j\leq n} }|S_{j}|\Vert_{p}$ to the
vector $(\Vert S_{j}\Vert_{2})_{1\leq j\leq n}$ we have\vspace*{1pt} to control $\sum
_{k=1}^{n}k^{-(p+1)/p}\|{\mathbf{E}}_{0}(S_{k})\|_{p}$ and $\sum_{k=1}%
^{n}k^{-(p+2\delta)/p}\Vert{\mathbf{E}}_{0}(S_{k}^{2})-{\mathbf
{E}}(S_{k}%
^{2})\Vert_{p/2}^{\delta}$.

In Section~\ref{martingalecase}, we study the case of stationary martingale
difference sequences showing that for all even powers $p\geq4$, inequality
(\ref{inedirect}) holds with \mbox{$\delta=2/(p-2)$}. This result is possible for
stationary martingale differences with the help of a special symmetrization
for martingales initiated by \citet{KwaWoy91}. In addition,
by using martingale approximation techniques, we obtain, for any even integer
$p$, another form of Rosenthal-type inequality, other than (\ref
{inedirect}), for
stationary adapted processes (see Section~\ref{sectionmartappro}), that
gives, for instance, better results for functionals of linear processes with
independent innovations.

We also investigate the situation when the conditional expectation with
respect to both the past and the future of the process is used. For instance,
when $p\geq4$ is an even integer, and the process is reversible, then
inequality (\ref{inedirect}) holds (see Theorem~\ref{stateven} and Corollary
\ref{Rev}) with $\delta=1$.

In Section~\ref{sectionBurk} we show that our inequalities imply the
Burkholder-type inequality as stated in Theorem 1 of
\citet{PelUteWu07}. For the sake of applications in Section \ref
{sectionindsum} we express
the terms that appear in our Rosenthal inequalities in terms of
individual summands.

Our paper is organized as follows. In Section~\ref{sectionmaxine}, we
prove a
new maximal inequality allowing us to relate the moments of the maximum of
partial sums of an adapted sequence, that is not necessarily
stationary, to
the moments of its partial sums. Maximal inequality (\ref{max}), combined
with moment estimates, allows us to obtain the Rosenthal-type inequalities
stated in Theorems~\ref{directprop} and~\ref{stateven} of Section
\ref{sectiongeneralpowers}. Section~\ref{sectionMA} is devoted to
Rosenthal-type inequalities for even powers for the special case of stationary
martingale differences and to an application to stationary processes
via a
martingale approximation technique. In Section \ref
{sectionappliexamples}, we
give other applications of the maximal inequalities stated in Section
\ref{sectionmaxine} and provide examples for which we compute the quantities
involved in the Rosenthal-type inequalities of Section~\ref{sectionmomentine}.
One of the applications presented in this section is a Bernstein inequality
for the maximum of partial sums for strongly mixing sequences that extends
the inequality in \citet{MerPelRio09}. The applications are
given to ARCH models, to functions of linear processes and reversible Markov
chains. In Section~\ref{densitysection}, we apply the inequality
(\ref{inedirect}) to estimate the random term of the ${\mathbf{L}}^{p}%
$-integrated risk of kernel estimators of the unknown marginal density
of a
stationary sequence that is assumed to be $\beta$-mixing in the weak sense
(see Definition~\ref{beta}). Some technical results are postponed
to the
\hyperref[app]{Appendix}.

\section{Maximal inequalities for adapted sequences}
\label{sectionmaxine}

The next proposition is a generalization of the well-known Doob maximal
inequality for martingales to adapted sequences. It states that the
moment of
order $p$ of the maximum of the partial sums of an adapted process can be
compared to the corresponding moment of the partial sum plus a
correction term
which is zero for martingale differences sequences. The proof is based on
convexity and chaining arguments.
\begin{proposition}
\label{maxinequality}Let $p>1$ and $q=p/(p-1)$. Let $Y_{i}$, $1\leq
i\leq
2^{r}$, be real random variables in ${\mathbf{L}}^{p}$, where $r$ is a positive
integer. Assume that the random variables are adapted to an
increasing\vadjust{\goodbreak}
filtration $({\mathcal{F}}_{i})_{i}$. Let $S_{r}=Y_{1}+\cdots+Y_{r}$.
Then the
following inequality holds:
%
%
\begin{equation}\label{inequality}%
\Bigl\Vert{\max_{1\leq i\le
2^{r}}}|S_{i}|\Bigr\Vert_{p}\leq q\Vert
S_{2^{r}}\Vert_{p}+q\sum_{l=0}^{r-1}
\Biggl(\sum_{k=1}^{2^{r-l}-1}\bigl\|{\mathbf{E}
}\bigl(S_{(k+1)2^{l}}-S_{k2^{l}}|{\mathcal{F}}_{k2^{l}}\bigr)\bigr\|_{p}^{p}
\Biggr)^{1/p} .\hspace*{-32pt}
\end{equation}
\end{proposition}
\begin{corollary}
\label{cormaxstat} In the stationary case, we get that for any integer
$r\geq1$,
\[
\Bigl\Vert{\max_{1\leq i\leq2^{r}}}|S_{i}|\Bigr\Vert_{p}\leq q\Vert
S_{2^{r}}\Vert_{p}+q2^{r/p}\sum_{l=0}^{r-1}2^{-l/p}\|{\mathbf
{E}}(S_{2^{l}%
}|{\mathcal{F}}_{0})\|_{p}.
\]
\end{corollary}
\begin{remark}
\label{Rem1inegalite} The inequality in Corollary~\ref{cormaxstat} easily
implies that%
%
%
\begin{equation} \label{maxdyatic}%
\Bigl\Vert{\max_{1\leq i\leq n}}|S_{i}|\Bigr\Vert_{p}\leq2q\max_{1\leq
m\leq n}\Vert S_{m}\Vert_{p}+(2^{1/p}q)n^{1/p}\sum_{l=0}^{r-1}2^{-l/p}%
\|{\mathbf{E}}(S_{2^{l}}|{\mathcal{F}}_{0})\|_{p} \hspace*{-20pt}
\end{equation}
for any integer $n\in[2^{r-1},2^{r}[$, where $r$ is a positive integer.
Moreover, due to the subadditivity of the sequence $ (\|{\mathbf{E}}%
(S_{n}|{\mathcal{F}}_{0})\|_{p})_{n\geq1}$, according to Lemma~\ref{lmasubadd},
we also have that for any positive integer $n$,
%
%
\begin{eqnarray}\label{cons1sub}%
\Bigl\Vert{\max_{1\leq i\leq n}}|S_{i}|\Bigr\Vert_{p}&\leq&2q\max_{1\leq
m\leq n}\Vert S_{m}\Vert_{p}\nonumber\\[-8pt]\\[-8pt]
&&{}+\biggl(q\frac{2^{2+2/p}}{2^{1+1/p}-1}\biggr)n^{1/p}\sum
_{j=1}^{n}j^{-1-1/p}\|{\mathbf{E}}(S_{j}|{\mathcal{F}}_{0})\|_{p}
.\nonumber
\end{eqnarray}
Inequalities (\ref{maxdyatic}) and (\ref{cons1sub}) are true, even if the
variables are not centered.
\end{remark}

For several applications involving exponential bounds we point out the
following proposition.
\begin{proposition}
\label{propmaxineproba} Let $p>1$ and $q=p/(p-1)$. Let $(Y_{i})_{i\geq
1}$ be
real random variables in ${\mathbf{L}}^{p}$. Assume that the random variables
are adapted to an increasing filtration $({\mathcal{F}}_{i})_{i}$. Let
$S_{n}=Y_{1}+\cdots+Y_{n}$. Let $\varphi$ be a nondecreasing, nonnegative,
convex and even function. Then for any positive real $x$, and any positive
integer $r$, the following inequality holds:
%
%
\begin{eqnarray}\label{inegalitemax}%
&&{\mathbf{P}}\Bigl({\max_{1\leq i\leq2^{r}}}|S_{i}|\geq2x\Bigr)\nonumber\\
&&\qquad\leq\frac{1}{\varphi
(x)}{\mathbf{E}}(\varphi(S_{2^{r}}))\\
&&\qquad\quad{}+q^{p}x^{-p}\Biggl(\sum_{l=0}^{r-1}%
\Biggl(\sum_{k=1}^{2^{r-l}-1}\bigl\|{\mathbf{E}}\bigl(S_{(k+1)2^{l}}-S_{k2^{l}%
}|{\mathcal{F}}_{k2^{l}}\bigr)\bigr\|_{p}^{p} \Biggr)^{1/p} \Biggr)^{p}
.\nonumber
\end{eqnarray}
Assume, in addition, that there exists a positive real $M$ such that
${\sup_{i}}\Vert Y_{i}\Vert_{\infty}\leq M$. Then for any positive real $x$,
and any
positive integer $r$, the following inequality holds:
%
%
\begin{eqnarray}\label{inegalitemax2proba}%
&&{\mathbf{P}}\Bigl({\max_{1\leq i\leq2^{r}}}|S_{i}|\geq4x\Bigr)\nonumber\\
&&\qquad\leq\frac{1}{\varphi
(x)}{\mathbf{E}}(\varphi(S_{2^{r}}))\\
&&\qquad\quad{}+q^{p}x^{-p}\Biggl(\sum_{l=0}^{r-1}%
\Biggl(\sum_{k=1}^{2^{r-l}-1}\bigl\|{\mathbf{E}}\bigl(S_{v+(k+1)2^{l}}-S_{v+k2^{l}%
}|{\mathcal{F}}_{k2^{l}}\bigr)\bigr\|_{p}^{p} \Biggr)^{1/p} \Biggr)^{p}
\nonumber
\end{eqnarray}
with $v=[x/M]$ (where \mbox{$[\cdot]$} denotes the integer part).
\end{proposition}
\begin{pf*}{Proofs of Propositions~\ref{maxinequality} and
\ref{propmaxineproba}}
Denote $S_{2^{r}}^{\ast}={\max_{1\leq i\leq2^{r}}}%
|S_{i}|$. For any $m\in[0,2^{r}-1]$, we have that
\[
S_{2^{r}-m}={\mathbf{E}}(S_{2^{r}}|{\mathcal{F}}_{2^{r}-m})-{\mathbf{E}%
}(S_{2^{r}}-S_{2^{r}-m}|{\mathcal{F}}_{2^{r}-m}) .
\]
So%
%
%
\begin{equation}\label{dec1max}\qquad
S_{2^{r}}^{\ast}\leq{\max_{0\leq m\leq2^{r}-1}}|{\mathbf{E}}(S_{2^{r}%
}|{\mathcal{F}}_{2^{r}-m})|+{\max_{0\leq m\leq2^{r}-1}}\bigl|{\mathbf
{E}}(S_{2^{r}%
}-S_{2^{r}-m}|{\mathcal{F}}_{2^{r}-m})\bigr|.
\end{equation}
Since $({\mathbf{E}}(S_{2^{r}}|{\mathcal{F}}_{k}))_{k\geq1}$ is a martingale,
we shall use Doob's maximal inequality to deal with the first term in the
right-hand side of (\ref{dec1max}). Hence, since $\varphi$ is a nondecreasing,
nonnegative, convex and even function, we get that
%
%
\begin{equation}\label{b1max}%
{\mathbf{P}} \Bigl({\max_{0\leq m\leq2^{r}-1}}|{\mathbf{E}}(S_{2^{r}%
}|{\mathcal{F}}_{2^{r}-m})|\geq x \Bigr)\leq\frac{1}{\varphi(x)}{\mathbf{E}%
}(\varphi(S_{2^{r}}))
\end{equation}
and also that
%
%
\begin{equation}\label{momentb1max}%
\Bigl\Vert{\max_{0\leq m\leq2^{r}-1}}|{\mathbf{E}}(S_{2^{r}}|{\mathcal
{F}}_{2^{r}%
-m})|\Bigr\Vert_{p}\leq q\Vert S_{2^{r}}\Vert_{p} .
\end{equation}
Write now $m$ in basis $2$ as follows:
\[
m=\sum_{i=0}^{r-1}b_{i}(m)2^{i}\qquad
\mbox{where $b_{i}(m)=0$ or
$b_{i}(m)=1$}.
\]
Set $m_{l}=\sum_{i=l}^{r-1}b_{i}(m)2^{i}$ and notice that for any $p\geq
1$, we
have
\[
\bigl|{\mathbf{E}}(S_{2^{r}}-S_{2^{r}-m}|{\mathcal{F}}_{2^{r}-m})\bigr|^{p}%
\leq\Biggl(\sum_{l=0}^{r-1}\bigl|{\mathbf{E}}(S_{2^{r}-m_{l+1}}-S_{2^{r}-m_{l}%
}|{\mathcal{F}}_{2^{r}-m})\bigr| \Biggr)^{p}.
\]
Hence setting
\[
\alpha_{l}= \Biggl(\sum_{k=1}^{2^{r-l}-1}\bigl\|{\mathbf{E}}\bigl(S_{(k+1)2^{l}}%
-S_{k2^{l}}|{\mathcal{F}}_{k2^{l}}\bigr)\bigr\|_{p}^{p} \Biggr)^{1/p}
\quad\mbox{and}\quad
\lambda_{l}=\frac{\alpha_{l}}{\sum_{l=0}^{r-1}\alpha_{l}} ,
\]
we get by convexity
\begin{eqnarray*}
&&\bigl|{\mathbf{E}}(S_{2^{r}}-S_{2^{r}-m}|{\mathcal{F}}_{2^{r}-m})\bigr|^{p}\\
&&\qquad\leq
\sum_{l=0}^{r-1}\lambda_{l}^{1-p}\bigl|{\mathbf{E}}(S_{2^{r}-m_{l+1}}%
-S_{2^{r}-m_{l}}|{\mathcal{F}}_{2^{r}-m})\bigr|^{p}\\
&&\qquad\leq \sum_{l=0}^{r-1}\lambda_{l}^{1-p}\bigl|{\mathbf{E}}\bigl(\bigl|{\mathbf{E}}%
(S_{2^{r}-m_{l+1}}-S_{2^{r}-m_{l}}|{\mathcal{F}}_{2^{r}-m_{l}}%
)\bigr| |{\mathcal{F}}_{2^{r}-m}\bigr)\bigr|^{p} .
\end{eqnarray*}
Now $m_{l}\neq m_{l+1}$, only if $b_{l}(m)=1$, and in that case $m_{l}%
=k_{m}2^{l}$ with $k_{m}$ odd. It follows that
\begin{eqnarray*}
&&\bigl|{\mathbf{E}}(S_{2^{r}-m_{l+1}}-S_{2^{r}-m_{l}}|
{\mathcal{F}}_{2^{r}-m_{l}})\bigr|\\
&&\qquad\leq
\max_{1\leq k\leq2^{r-l},k\ \mathrm{odd}}\bigl|{\mathbf{E}}\bigl(S_{2^{r}%
-(k-1)2^{l}}-S_{2^{r}-k2^{l}}|{\mathcal{F}}_{2^{r}-k2^{l}}\bigr)\bigr|\\
&&\qquad:= A_{r,l} .
\end{eqnarray*}
Hence, using the fact that if $|X|\leq|Y|$, then ${\mathbf
{E}}(|X||{\mathcal{F}%
})\leq{\mathbf{E}}(|Y||{\mathcal{F}})$, we get that
\[
\Bigl\Vert\max_{0\leq m\leq2^{r}-1}\bigl|{\mathbf{E}}(S_{2^{r}}-S_{2^{r}-m}%
|{\mathcal{F}}_{2^{r}-m})\bigr|\Bigr\Vert_{p}^{p}
\leq\sum_{l=0}^{r-1}\lambda_{l}%
^{1-p}{\mathbf{E}}\Bigl({\max_{0\leq m\leq2^{r}-1}}|{\mathbf{E}}(A_{r,l}%
|{\mathcal{F}}_{2^{r}-m})|\Bigr)^{p} .
\]
Notice that $({\mathbf{E}}(A_{r,l}|{\mathcal{F}}_{k}))_{k\geq1}$ is a
martingale and by Doob's maximal inequality, we obtain
\[
{\mathbf{E}}\Bigl(\max_{0\leq m\leq2^{r}-1}{\mathbf{E}}(A_{r,l}|{\mathcal{F}%
}_{2^{r}-m})\Bigr)^{p}\leq q^{p}\Vert A_{r,l}\Vert_{p}^{p}\leq q^{p}\alpha
_{l}%
^{p} .
\]
Using the definition of $\lambda_{l}$, it follows that
%
%
\begin{eqnarray}\label{b2maxmomentproba}
&&\Bigl\Vert\max_{0\leq m\leq2^{r}-1}\bigl|{\mathbf{E}}(S_{2^{r}}-S_{2^{r}-m}%
|{\mathcal{F}}_{2^{r}-m})\bigr|\Bigr\Vert_{p}^{p}\nonumber\\
&&\qquad\leq q^{p} \Biggl(\sum_{l=0}^{r-1}%
\alpha_{l} \Biggr)^{p}\\
&&\qquad\leq q^{p}\Biggl(\sum_{l=0}^{r-1}\Biggl(\sum_{k=1}^{2^{r-l}-1}\bigl\|{\mathbf{E}%
}\bigl(S_{(k+1)2^{l}}-S_{k2^{l}}|{\mathcal{F}}_{k2^{l}}\bigr)
\bigr\|_{p}^{p} \Biggr)^{1/p}%
\Biggr)^{p} .\nonumber
\end{eqnarray}

Starting from (\ref{dec1max}) and using (\ref{b1max}) and
(\ref{b2maxmomentproba}) combined with Markov's inequality, inequality
(\ref{inegalitemax}) of Proposition~\ref{propmaxineproba} follows. To
end the
proof of Proposition~\ref{maxinequality}, we start from (\ref{dec1max}) and
consider bounds (\ref{momentb1max}) and (\ref{b2maxmomentproba}).

We turn now to the proof of inequality (\ref{inegalitemax2proba}). We start
from (\ref{dec1max}) and we write that
\begin{eqnarray*}
S_{2^{r}}^{\ast} & \leq &{\max_{0\leq m\leq2^{r}-1}}|{\mathbf{E}}(S_{2^{r}%
}|{\mathcal{F}}_{2^{r}-m})|\\
&&{}+\max_{0\leq m\leq2^{r}-1}\bigl|{\mathbf
{E}}(S_{2^{r}%
+v}-S_{2^{r}+v-m}|{\mathcal{F}}_{2^{r}-m})\bigr|\\
&&{} +\max_{0\leq m\leq2^{r}-1}\bigl|{\mathbf
{E}}(S_{2^{r}+v}-S_{2^{r}}|{\mathcal{F}%
}_{2^{r}-m})\bigr|\\
&&{}+{\max_{0\leq m\leq2^{r}-1}}\bigl|{\mathbf
{E}}(S_{2^{r}+v-m}-S_{2^{r}%
-m}|{\mathcal{F}}_{2^{r}-m})\bigr| .
\end{eqnarray*}
By the fact that the variables are uniformly bounded by $M$, we then derive
that
\[
S_{2^{r}}^{\ast}\leq{\max_{0\leq m\leq2^{r}-1}}|{\mathbf{E}}(S_{2^{r}%
}|{\mathcal{F}}_{2^{r}-m})|+\max_{0\leq m\leq2^{r}-1}\bigl|{\mathbf
{E}}(S_{2^{r}%
+v}-S_{2^{r}+v-m}|{\mathcal{F}}_{2^{r}-m})\bigr|+2vM .
\]
Since $vM\leq x$, it follows that
%
%
\begin{eqnarray}\label{b3maxmomentproba}
{\mathbf{P}}(S_{2^{r}}^{\ast}\geq4x)&\leq&{\mathbf{P}}\Bigl({\max_{0\leq m\leq
2^{r}%
-1}}|{\mathbf{E}}(S_{2^{r}}|{\mathcal{F}}_{2^{r}-m})|\geq
x\Bigr)\nonumber\\[-8pt]\\[-8pt]
&&{}+{\mathbf{P}}\Bigl(\max_{0\leq m\leq2^{r}-1}\bigl|{\mathbf
{E}}(S_{2^{r}+v}-S_{2^{r}%
+v-m}|{\mathcal{F}}_{2^{r}-m})\bigr|\geq x\Bigr) .\nonumber
\end{eqnarray}
Using chaining arguments, convexity and Doob's maximal inequality, as
above, we infer that for any $p>1$,
%
%
\begin{eqnarray}\label{b4maxmomentproba}
&& \Bigl\Vert\max_{0\leq m\leq2^{r}-1}\bigl|{\mathbf{E}}(S_{2^{r}+v}-S_{2^{r}%
+v-m}|{\mathcal{F}}_{2^{r}-m})\bigr|\Bigr\Vert_{p}^{p}\nonumber\\[-8pt]\\[-8pt]
&&\qquad \leq q^{p}\Biggl(\sum_{l=0}^{r-1} \Biggl(\sum_{k=1}^{2^{r-l}-1}
\bigl\|{\mathbf{E}%
}\bigl(S_{v+(k+1)2^{l}}-S_{v+k2^{l}}|{\mathcal{F}}_{k2^{l}}\bigr)\bigr\|_{p}^{p}%
\Biggr)^{1/p} \Biggr)^{p} .\nonumber
\end{eqnarray}
Starting from (\ref{b3maxmomentproba}) and using (\ref{b1max}) and
(\ref{b4maxmomentproba}) combined with Markov's inequality, inequality
(\ref{inegalitemax2proba}) of Proposition~\ref{propmaxineproba} follows.
\end{pf*}

\section{Moment inequalities for the maximum of partial sums under projective
conditions}
\label{sectionmomentine}

\subsection{Rosenthal-type inequalities for stationary processes}
\label{sectiongeneralpowers}

Using a direct approach that combines maximal inequality (\ref{cons1sub})
and Lemma~\ref{cross}, we obtain the following Rosenthal inequality for
the maximum of the partial sums of a stationary process for all powers $p>2$.
\begin{theorem}
\label{directprop} Let $p>2$ be a real number, and let $(X_{i})_{i\in
\mathbb{Z}}$ be an adapted stationary sequence in the sense of Notation
\ref{notation1}.\vadjust{\goodbreak} Then, for any positive integer~$n$, the following inequality
holds:
\begin{eqnarray*}
{\mathbf{E}} \Bigl({\max_{1\leq j\leq n}}|S_{j}|^{p} \Bigr)&\ll&
n{\mathbf{E}%
}(|X_{1}|)^{p}+cn\Biggl( \sum_{k=1}^{n}\frac{1}{k^{1+1/p}}\|{\mathbf{E}}%
_{0}(S_{k})\|_{p}\Biggr)^{p}\\
&&{}+n\Biggl( \sum_{k=1}^{n}\frac{1}{k^{1+2\delta
/p}}\Vert{\mathbf{E}}_{0}(S_{k}^{2})\Vert_{p/2}^{\delta}\Biggr)
^{p/(2\delta)},
\end{eqnarray*}
where $\delta=\min(1,1/(p-2))$ and $c=1$. When $p\geq4$ we can take
$c=0$ by
enlarging the constant involved.
\end{theorem}

With applications to Markov processes in mind, by conditioning with
respect to
both the future and the past of the process, our next result gives an
alternative inequality to the one given in Theorem~\ref{directprop}
when $p$
is an even integer. For this case, the power $\delta$ appearing in Theorem
\ref{directprop} is always equal to one. Before stating the result, we first
introduce the following notation to define the additional nonincreasing
filtration that we consider.
\begin{notation}
\label{notation2} Let\vspace*{2pt} $\mathcal{\bar{F}}_{0}$ be a $\sigma$-algebra of
$\mathcal{A}$ satisfying $T^{-1}(\mathcal{\bar{F}}_{0})\subseteq
\mathcal{\bar{F}}_{0}$. We then define the nonincreasing filtration
$(\mathcal{\bar{F}}_{i})_{i\in\mathbb{Z}}$ by $\mathcal{\bar{F}}_{i}%
=T^{-i}(\mathcal{\bar{F}}_{0})$. In what follows, we use the notation
${\bar{\mathbf{E}}}_{k}(Y)={\mathbf{E}}(Y|\mathcal{\bar{F}}_{k})$.
\end{notation}
\begin{theorem}
\label{stateven} Let $p\geq4$ be an even integer, and let $X_{0}$ be a real
valued random variable such that $\Vert X_{0}\Vert_{p}<\infty$ and measurable
with respect to ${\mathcal{F}}_{0}$ and to ${\mathcal{\bar{F}}}_{0}$. We
construct the stationary sequence $(X_{i})_{i\in\mathbb{Z}}$ as in Notation
\ref{notation1}. Then for any integer~$n$,
\begin{eqnarray*}
{\mathbf{E}} \Bigl({\max_{1\leq j\leq n}}|S_{j}|^{p} \Bigr) &\ll& n{\mathbf{E}%
}(|X_{1}|^{p})+n \Biggl(\sum_{k=1}^{n}\frac{1}{k^{1+1/p}} \bigl(\Vert{\mathbf{E}%
}_{0}(S_{k})\Vert_{p}+\Vert{\bar{\mathbf{E}}}_{k+1}(S_{k})\Vert_{p}%
\bigr) \Biggr)^{p}\\
&&{} +n \Biggl(\sum_{k=1}^{n}\frac{1}{k^{1+2/p}} \bigl(\Vert{\mathbf{E}}_{0}%
(S_{k}^{2})\Vert_{p/2}+\Vert{\bar{\mathbf{E}}}_{k+1}(S_{k}^{2})\Vert
_{p/2} \bigr) \Biggr)^{p/2} .
\end{eqnarray*}
\end{theorem}

As a corollary to the proof of Theorem~\ref{stateven}, we obtain:
\begin{theorem}
\label{pinteger}Let $p\geq4$ be a real number, and let $X_{0}$ be a real-valued
random variable such that $\Vert X_{0}\Vert_{p}<\infty$ and measurable with
respect to ${\mathcal{F}}_{0}$ and to~${\mathcal{\bar{F}}}_{0}$. We construct
the stationary sequence $(X_{i})_{i\in\mathbb{Z}}$ as in Notation
\ref{notation1}. Then for any integer~$n$,
\[
{\mathbf{E}} \Bigl({\max_{1\leq j\leq n}}|S_{j}|^{p} \Bigr)\ll n{\mathbf{E}%
}(|X_{1}|^{p})+n \Biggl(\sum_{k=1}^{n}\frac{1}{k^{1+1/p}} \bigl(\Vert{\mathbf{E}%
}_{0}(S_{k}^{2})\Vert_{p/2}^{1/2}+\Vert{\bar{\mathbf{E}}}_{k+1}(S_{k}%
^{2})\Vert_{p/2}^{1/2} \bigr) \Biggr)^{p} .
\]
\end{theorem}

This theorem is also valid for $2<p<4$. In this range, however,
according to
item (2) of Comment~\ref{commentgreaterdelta}, Theorem~\ref{directprop} gives
better bounds.
\begin{pf*}{Proof of Theorem~\ref{directprop}}
The proof of this theorem
is based on dyadic induction and involves several steps. With the notation
$a_{n}=\|S_{n}\|_{p}$ we shall establish a recurrence formula: for any
positive integer $r$,%
%
%
\begin{equation} \label{recurrence}%
a_{2n}^{p}\leq2a_{n}^{p}+2c_{1}a_{n}^{p-1}\Vert{\mathbf{E}}_{0}(S_{n}%
)\Vert_{p}+2c_{2}a_{n}^{p-2\delta}\Vert{\mathbf{E}}_{0}(S_{n}^{2})\Vert
_{p/2}^{\delta} ,
\end{equation}
where $c_{1}$ and $c_{2}$ are positive constants depending only on $p$. Before
proving it, let us show that (\ref{recurrence}) implies our result.
\begin{lemma}
\label{reclemma}Assume that for some $0<\delta\leq1$ the recurrence formula
(\ref{recurrence}) holds. Then inequalities (\ref{dyaticform}) and
(\ref{inedirect}) hold with the same $\delta$.
\end{lemma}

Let us prove the lemma. From inequality (\ref{recurrence}), by
recurrence on
the first term, we obtain for any positive integer $r$,
%
%
\begin{eqnarray}\label{recurrence2}%
a_{2^{r}}^{p}&\leq&2^{r} \Biggl(a_{2^{0}}^{p}+c_{1}\sum
_{k=0}^{r-1}2^{-k}a_{2^{k}%
}^{p-1}\Vert{\mathbf{E}}_{0}(S_{2^{k}})\Vert_{p}\nonumber\\[-8pt]\\[-8pt]
&&\hspace*{2.2pt}\hphantom{2^{r} \Biggl(}
{}+c_{2}\sum_{k=0}^{r-1}%
2^{-k}a_{2^{k}}^{p-2\delta}\Vert{\mathbf{E}}_{0}(S_{2^{k}}^{2})\Vert
_{p/2}^{\delta} \Biggr).\nonumber
\end{eqnarray}

We shall establish first inequality (\ref{dyaticform}). Due to maximal
inequality (\ref{maxdyatic}), it suffices to prove that the inequality is
satisfied for $\max_{1\leq j\leq n}{\mathbf{E}}(|S_{j}|^{p})$ instead of
${\mathbf{E}} (\max_{1\leq j\leq n}|S_{j}|^{p} )$.

The proof is divided in several steps. The goal is to establish that
for any
positive integer $r$ and any integer $n$ such that $2^{r-1}\leq n<2^{r}$,
%
%
\begin{eqnarray}\label{ineqtoshow}%
\max_{1\leq j\leq n}{\mathbf{E}} (|S_{j}|^{p} )&\ll& n{\mathbf{E}%
}(|X_{1}|)^{p}+cn \Biggl(\sum_{k=0}^{r-1}2^{-k/p}\Vert{\mathbf
{E}}_{0}(S_{2^{k}%
})\Vert_{p} \Biggr)^{p}\nonumber\\[-8pt]\\[-8pt]
&&{}+n \Biggl(\sum_{k=0}^{r-1}2^{-2k\delta/p}\Vert{\mathbf{E}%
}_{0}(S_{2^{k}}^{2})\Vert_{p/2}^{\delta} \Biggr)^{p/2\delta} .\nonumber
\end{eqnarray}
With the notation $ B_{r}=\max_{0\leq k\leq r}(a_{2^{k}}%
^{p}/2^{k})$, starting from (\ref{recurrence2}), we get
\begin{eqnarray*}
B_{r}&\leq& a_{2^{0}}^{p}+c_{1}B_{r}^{1-1/p}\sum_{k=0}^{r-1}2^{-k/p}%
\Vert{\mathbf{E}}_{0}(S_{2^{k}})\Vert_{p}\\
&&{}+c_{2}B_{r}^{1-2\delta/p}\sum
_{k=0}^{r-1}2^{-2k\delta/p}\Vert{\mathbf{E}}_{0}(S_{2^{k}}^{2})\Vert
_{p/2}^{\delta} .
\end{eqnarray*}
Therefore, taking into account that either $B_{r}\leq3a_{2^{0}}^{p}$ or
\[
B_{r}^{1/p}\leq3c_{1}\sum_{k=0}^{r-1}2^{-k/p}\Vert{\mathbf
{E}}_{0}(S_{2^{k}%
})\Vert_{p}\quad \mbox{or}\quad B_{r}^{2\delta/p}\leq3c_{2}\sum_{k=0}^{r-1}2^{-2k\delta
/p} \Vert{\mathbf{E}}_{0}(S_{2^{k}}^{2})\Vert_{p/2}^{\delta},
\]
we derive that
%
%
\begin{eqnarray}\label{boundpropdyatic2}%
a_{2^{r}}^{p}&\leq&2^{r} \Biggl(3a_{2^{0}}^{p}+ \Biggl(3c_{1}\sum_{k=0}%
^{r-1}2^{-k/p}\Vert{\mathbf{E}}_{0}(S_{2^{k}})\Vert_{p}
\Biggr)^{p}\nonumber\\[-8pt]\\[-8pt]
&&\hspace*{0pt}\hphantom{2^{r} \Biggl(}
{} + \Biggl(3c_{2}\sum_{k=0}^{r-1}2^{-2k\delta/p}\Vert{\mathbf{E}}_{0}(S_{2^{k}%
}^{2})\Vert_{p/2}^{\delta} \Biggr)^{p/2\delta} \Biggr) .\nonumber
\end{eqnarray}
Let now $2^{r-1}\leq n<2^{r}$, and write its binary expansion,
%
%
\begin{equation} \label{binaryexpansion}%
n=\sum_{k=0}^{r-1}2^{k}b_{k}\qquad\mbox{where }b_{r-1}=1\mbox{ and }b_{k}%
\in\{0,1\}\qquad\mbox{for }k=0,\ldots,r-2 .\hspace*{-28pt}
\end{equation}
Notice that
\begin{eqnarray*}
S_{n}&=&\sum_{k=0}^{r-1}b_{k}T_{2^{k}}\qquad\mbox{where }T_{2^{k}}=\sum
_{i=n_{k-1}%
+1}^{n_{k}}X_{i},\\
n_{k}&=&\sum_{j=0}^{k}b_{j}2^{j}\quad\mbox{and}\quad n_{-1}=0.%
\end{eqnarray*}
Hence, by stationarity,
\[
\Vert S_{n}\Vert_{p}\leq\sum_{k=0}^{r-1}b_{k}\Vert T_{2^{k}}\Vert
_{p}\leq
\sum_{k=0}^{r-1}b_{k}\Vert S_{2^{k}}\Vert_{p} .
\]
Then, by using (\ref{boundpropdyatic2}) and the fact that $\sum_{k=0}%
^{r-1}b_{k}2^{k/p}\leq2^{r/p}/(1-2^{-1/p})$, we derive inequality
(\ref{ineqtoshow}) for ${\mathbf{E}}(|S_{n}|^{p})$ and also for $\max
_{1\leq
j\leq n}{\mathbf{E}} (|S_{j}|^{p} )$. Inequality (\ref{dyaticform}%
) follows now by maximal inequality (\ref{maxdyatic}).

We indicate now how to derive from (\ref{dyaticform}) the inequality
stated in
Theorem~\ref{directprop}. Notice that, by stationarity, for any
integers $i$
and $j$,
\[
\Vert{\mathbf{E}}_{0}(S_{i+j})\Vert_{p}\leq\Vert{\mathbf{{E}}}_{0}(S_{i}
)\Vert_{p}+\Vert{\mathbf{{E}}}_{0}(S_{j})\Vert_{p} ,
\]
and also that for any $0<\delta\leq1$,
\[
\Vert{\mathbf{E}}_{0}(S_{i+j}^{2})\Vert_{p/2}^{\delta}\leq2^{\delta}%
\Vert{\mathbf{E}}_{0}(S_{i}^{2})\Vert_{p/2}^{\delta}+2^{\delta}\Vert
{\mathbf{E}}_{0}(S_{j}^{2})\Vert_{p/2}^{\delta} .
\]
Using item (1) of Lemma~\ref{lmasubadd}, it follows that
\[
\sum_{k=0}^{r-1}2^{-k/p}\Vert{\mathbf{E}}_{0}(S_{2^{k}})\Vert_{p}\ll\sum
_{k=1}^{n}k^{-1-1/p}\Vert{\mathbf{E}}_{0}(S_{k})\Vert_{p}
\]
and%
%
%
\begin{equation}\label{majdyadelta}%
\sum_{k=0}^{r-1}2^{-2k\delta/p}\Vert{\mathbf{E}}_{0}(S_{2^{k}}^{2})\Vert
_{p/2}^{\delta}\ll\sum_{k=1}^{n}k^{-1-2\delta/p}\Vert{\mathbf{E}}_{0}%
(S_{k}^{2})\Vert_{p/2}^{\delta} .
\end{equation}
The results follows by the above considerations via inequality
(\ref{dyaticform}).
\end{pf*}
\begin{pf*}{End of the proof of Theorem~\ref{directprop}} It remains to
establish the recurrence formula (\ref{recurrence}). We
divide\vspace*{1pt} the proof into three cases according to the values
of $p$. Denote $\bar{S}_{n}=X_{n+1}+\cdots+X_{2n}$.

The case $2<p\leq3$ was discussed in \citet{Rio09}. We give here
a shorter
alternative proof. We apply inequality (\ref{Rio}) of Lemma~\ref{cross} with
$x=S_{n}$ and $y=\bar{S}_{n}$. Then, by taking the expectation and using
stationarity and properties of conditional expectation, we obtain
\[
{\mathbf{E}}(|S_{2n}|)^{p}\leq2{\mathbf{E}}(|S_{n}|)^{p}+p{\mathbf{E}}%
(|S_{n}|^{p-1}\operatorname{sign}(S_{n}){\mathbf{E}}_{n}(\bar
{S}_{n}))+C_{p}%
^{2}{\mathbf{E}}(|S_{n}|^{p-2}{\mathbf{E}}_{n}(\bar{S}_{n}^{2})) ,
\]
where $C_{p}^{2}=p(p-1)/2$. This inequality combined with H\"{o}lder's
inequality gives
\[
a_{2n}^{p}\leq2a_{n}^{p}+pa_{n}^{p-1}\Vert{\mathbf{E}}_{0}(S_{n})\Vert
_{p}+C_{p}^{2}a_{n}^{p(1-2/p)}\Vert{\mathbf{E}}_{0}(S_{n}^{2})\Vert
_{p/2} ,
\]
and therefore (\ref{recurrence}) holds with $\delta=1$, $c_{1}=2^{-1}p$ and
$c_{2}=2^{-1}C_{p}^{2}$.

Assume now that $p\in\ ]3,4[$. Using inequality (\ref{Rio34}) of Lemma
\ref{cross} with $x=S_{n}$ and $y=\bar{S}_{n}$, taking the expectation, and
using Lemma~\ref{basic}, we get by stationarity that for any positive integer
$r$,
\begin{eqnarray*}
a_{2n}^{p}&\leq&2a_{n}^{p}+pa_{n}^{p-1}\Vert{\mathbf{E}}_{0}(S_{n})\Vert
_{p}+C_{p}^{2}a_{n}^{p(1-2/p)}\Vert{\mathbf{E}}_{0}(S_{n}^{2})\Vert
_{p/2}\\
&&{}+2p(p-2)^{-1}a_{n}^{p-2/(p-2)}\|{\mathbf{E}}_{0}(S_{n}^{2}%
)\|_{p/2}^{1/(p-2)} .
\end{eqnarray*}
Since $\Vert{\mathbf{E}}_{0}(S_{n}^{2})\Vert_{p/2}\leq
a_{n}^{2(p-3)/(p-2)}%
\Vert{\mathbf{E}}_{0}(S_{n}^{2})\Vert_{p/2}^{1/(p-2)}$, it follows that
\[
a_{2n}^{p}\leq2a_{n}^{p}+pa_{n}^{p-1}\Vert{\mathbf{E}}%
_{0}(S_{n})\Vert_{p}+4pa_{n}^{p-2/(p-2)}\|{\mathbf{E}}_{0}(S_{n}^{2}%
)\|_{p/2}^{1/(p-2)} .
\]
It follows that (\ref{recurrence}) holds with $\delta=1/(p-2)$, $c_{1}%
=2^{-1}p$ and $c_{2}=2p$.\vadjust{\goodbreak}

It remains to prove inequality (\ref{inedirect}) for $p\geq4$. Using
inequality (\ref{int}) of Lemma~\ref{cross} with $x=S_{n}$ and $y=\bar
{S}_{n}%
$, and taking the expectation, we get by stationarity that
\[
a_{2n}^{p}\leq2a_{n}^{p}+4^{p}{\mathbf{E}} (|S_{n}|^{p-1}|\bar{S}%
_{n}|+|\bar{S}_{n}|^{p-1}|S_{n}| ) .
\]
Using Lemma~\ref{basic} together with stationarity, it follows that
\[
{\mathbf{E}}(|S_{n}||\bar{S}_{n}|^{p-1})\leq a_{n}^{p-2/(p-2)}\|{\mathbf
{E}%
}_{0}(S_{n}^{2})\|_{p/2}^{1/(p-2)}
\]
and that
\[
{\mathbf{E}}(|S_{n}|^{p-1}|\bar{S}_{n}|)\leq a_{n}^{p-1}\|{\mathbf{E}}%
_{0}(S_{n}^{2})\|_{p/2}^{1/2}\leq a_{n}^{p-2/(p-2)}\|{\mathbf{E}}_{0}%
(S_{n}^{2})\|_{p/2}^{1/(p-2)} .
\]
From these estimates we deduce
\[
a_{2n}^{p}\leq2a_{n}^{p}+2(4^{p})a_{n}^{p-2/(p-2)}\Vert{\mathbf{E}}_{0}%
(S_{n}^{2})\Vert_{p/2}^{1/(p-2)} ,
\]
and then (\ref{recurrence}) holds with $\delta=1/(p-2)$, $c_{1}=0$ and
$c_{2}=4^{p}$. Therefore, in this case (\ref{ineqtoshow}) holds with $c=0$.
Then by maximal inequality (\ref{maxdyatic}), inequality
(\ref{dyaticform}) holds with $c=1$. We show now that, in this case, the
second term in the inequality can be bounded up to a multiplicative constant
by the third term. By Jensen's inequality and since in this case
$\delta<1/2$,
we have
%
%
\begin{eqnarray} \label{comp}%
\sum_{k=0}^{r-1}2^{-k/p}\|{\mathbf{E}}_{0}(S_{2^{k}})\|_{p}&\leq&\Biggl(\sum
_{k=0}^{r-1}2^{-(2k/p)(1/2)}\|{\mathbf{E}}_{0}(S_{2^{k}}^{2})\|_{p/2}%
^{1/2} \Biggr)^{2}\nonumber\\[-8pt]\\[-8pt]
&\leq&\Biggl(\sum_{k=0}^{r-1}2^{-(2k/p)\delta}\Vert{\mathbf{E}%
}_{0}(S_{2^{k}}^{2})\Vert_{p/2}^{\delta} \Biggr)^{1/\delta} .\nonumber
\end{eqnarray}
We then finish the proof by using (\ref{majdyadelta}).
\end{pf*}
\begin{pf*}{Proof of Theorem~\ref{stateven}}
Denote $\bar{S}%
_{n}=X_{n+1}+\cdots+X_{2n}$. Starting from inequality (\ref{evenint}) of
Lemma~\ref{cross} applied with $x=S_{n}$ and $y=\bar{S}_{n}$ and using the
notation $a_{n}=\|S_{n}\|_{p}$, by stationarity, we get that
%
%
\begin{eqnarray}\label{b1a2n}%
a_{2n}^{p}&\leq&2a_{n}^{p}+p \bigl({\mathbf{E}}(S_{n}^{p-1}{\bar{S}}%
_{n})+{\mathbf{E}}(S_{n}{\bar{S}}_{n}^{p-1}) \bigr)\nonumber\\[-8pt]\\[-8pt]
&&{}+2^{p} \bigl({\mathbf{E}%
}(S_{n}^{p-2}\bar{S}_{n}^{2})+{\mathbf{E}}(\bar{S}_{n}^{p-2}S_{n}%
^{2}) \bigr) .\nonumber
\end{eqnarray}
By using H\"{o}lder's inequality and recurrence, we then derive that
for any
positive integer $r$,
\begin{eqnarray*}
a_{2^{r}}^{p} & \leq & 2^{r}a_{2^{0}}^{p}+2^{-1}p\sum_{k=0}^{r-1}%
2^{r-k}a_{2^{k}}^{p-1} \bigl(\Vert{\mathbf{E}}_{0}(S_{2^{k}})\Vert_{p}%
+\Vert{\bar{\mathbf{E}}}_{2^{k}+1}(S_{2^{k}})\Vert_{p} \bigr)\\
&&{} +2^{p-1}\sum_{k=0}^{r-1}2^{r-k}a_{2^{k}}^{p-2} \bigl(\Vert{\mathbf{E}}%
_{0}(S_{2^{k}}^{2})\Vert_{p/2}+\Vert{\bar{\mathbf
{E}}}_{2^{k}+1}(S_{2^{k}}%
^{2})\Vert_{p/2} \bigr) .
\end{eqnarray*}
By using the arguments of the proof of Lemma~\ref{reclemma}, we get for
$2^{r-1}\leq n<2^{r}$,
\begin{eqnarray*}
{\mathbf{E}} \Bigl({\max_{1\leq j\leq n}}|S_{j}|^{p} \Bigr) & \ll & n{\mathbf{E}%
}(|X_{1}|^{p})+n \Biggl(\sum_{k=1}^{r-1}2^{-k/p} \bigl(\Vert{\mathbf{E}}%
_{0}(S_{2^{k}})\Vert_{p}+\Vert{\bar{\mathbf{E}}}_{k+1}(S_{2^{k}})\Vert
_{p} \bigr) \Biggr)^{p}\\
&&{} +n \Biggl(\sum_{k=1}^{n}2^{-2k/p} \bigl(\Vert{\mathbf{E}}_{0}(S_{2^{k}}%
^{2})\Vert_{p/2}+\Vert{\bar{\mathbf{E}}}_{k+1}(S_{2^{k}}^{2})\Vert
_{p/2} \bigr) \Biggr)^{p/2} .
\end{eqnarray*}
Noticing, in addition, that, by stationarity, for any integer $i$ and $j$,
%
%
\begin{eqnarray}\label{*}%
\Vert{\mathbf{E}}_{0}(S_{i+j})\Vert_{p}&\leq&\Vert{\mathbf{E}}_{0}(S_{i}%
)\Vert_{p}+\Vert{\mathbf{E}}_{0}(S_{j})\Vert_{p},\nonumber\\
\Vert{\bar
{\mathbf{E}}}_{i+j+1}(S_{i+j})\Vert_{p}&\leq&\Vert{\bar{\mathbf{E}}}_{i+1}(S_{i}%
)\Vert_{p}+\Vert{\bar{\mathbf{E}}}_{j+1}(S_{j})\Vert_{p} ,
\\
\Vert{\bar{\mathbf{E}}}_{i+j+1}(S_{i+j}^{2})\Vert_{p/2}&\leq&2\Vert
{\bar{\mathbf{E}}}_{i+1}(S_{i}^{2})\Vert_{p/2}+2\Vert{\bar{\mathbf{E}}}%
_{j+1}(S_{j}^{2})\Vert_{p/2} \nonumber
\end{eqnarray}
and that
%
%
\begin{equation}\label{**}%
\Vert{\mathbf{E}}_{0}(S_{i+j}^{2})\Vert_{p/2}\leq2\Vert{\mathbf{E}}_{0}%
(S_{i}^{2})\Vert_{p/2}+2\Vert{\mathbf{E}}_{0}(S_{j}^{2})\Vert_{p/2} .
\end{equation}
We obtain the desired result by using Lemma~\ref{lmasubadd}.
\end{pf*}
\begin{pf*}{Proof of Theorem~\ref{pinteger}} To prove this
theorem we apply inequality (\ref{int}) of Lemma~\ref{cross} with
$x=S_{n}$ and $y=\bar{S}_{n}$, where $\bar{S}_{n}=X_{n+1}+\cdots
+X_{2n}$. With
the notation $a_{n}=\|S_{n}\|_{p}$, we then have by stationarity that
\[
a_{2n}^{p}\leq2a_{n}^{p}+4^{p} \bigl({\mathbf{E}}(|S_{n}|^{p-1}|{\bar{S}}%
_{n}|)+{\mathbf{E}}(|S_{n}||{\bar{S}}_{n}|^{p-1}) \bigr) .
\]
By conditioning and then applying Jensen's inequality followed by H\"{o}lder's
inequality we obtain
\begin{eqnarray*}
a_{2n}^{p} & \leq & 2a_{n}^{p}+4^{p} \bigl({\mathbf{E}}(|S_{n}|^{p-1}%
{\mathbf{E}}_{n}^{1/2}({\bar{S}}_{n}^{2}))+{\mathbf{E}}(|{\bar{S}}_{n}%
|^{p-1}{\bar{\mathbf{E}}}_{n+1}^{1/2}(S_{n}^{2})) \bigr) \\
& \leq & 2a_{n}^{p}+4^{p} a_{n}^{p-1} \bigl(\Vert{\mathbf{E}}_{0}(S_{n}%
^{2})\Vert_{p/2}^{1/2}+\Vert{\bar{\mathbf{E}}}_{n+1}(S_{n}^{2})\Vert
_{p/2}^{1/2} \bigr) .
\end{eqnarray*}
By recurrence, we then derive that for any positive integer $r$,
\[
a_{2^{r}}^{p}\leq2^{r} \Biggl(a_{0}^{p}+2^{2p-1}\sum_{k=0}^{r-1}2^{-k}a_{2^{k}
}^{p-1} \bigl(\Vert{\mathbf{E}}_{0}(S_{2^{k}}^{2})\Vert_{p/2}^{1/2}%
+\Vert{\bar{\mathbf{E}}}_{2^{k}+1}(S_{2^{k}}^{2})\Vert_{p/2}^{1/2}%
\bigr) \Biggr) .
\]
The proof is completed by the arguments developed in the proof of Lemma~\ref{reclemma}
and by using Lemma~\ref{lmasubadd} via inequalities
(\ref{*}) and (\ref{**}) .
\end{pf*}

\subsection{Relation with the Burkholder-type inequality}
\label{sectionBurk}

The next lemma shows how to compare $\Vert{\mathbf
{E}}_{0}(S_{n}^{2})\Vert_{p/2}$
with quantities involving only $\Vert{\mathbf{E}}_{0}(S_{n})\Vert_{p}$.
\begin{lemma}
\label{comparison} Let $p\geq2$ be a real number, and let $(X_{n})$ be an
adapted stationary sequence in the sense of Notation~\ref{notation1}. Then,
for any positive integer~$n$,
%
%
\begin{equation}\label{comp2}%
\Vert{\mathbf{E}}_{0}(S_{n}^{2})\Vert_{p/2}\ll n\Vert{\mathbf{E}}_{0}%
(X_{1}^{2})\Vert_{p/2}+n\Biggl( \sum_{j=1}^{n}\frac{\Vert{\mathbf{E}}%
_{0}(S_{j})\Vert_{p}}{j^{3/2}}\Biggr) ^{2} .
\end{equation}
\end{lemma}

As a consequence of the above lemma, we get that for any $0<\delta\leq
1$ and
any real $p>2$,
\[
n\Biggl( \sum_{j=1}^{n}\frac{\Vert{\mathbf{E}}_{0}(S_{j}^{2})\Vert
_{p/2}^{\delta}}{j^{1+2\delta/p}}\Biggr) ^{p/(2\delta)}\ll n^{p/2}%
\Vert{\mathbf{E}}_{0}(X_{1}^{2})\Vert_{p/2}^{p/2}+n^{p/2}\Biggl( \sum
_{j=1}^{n}\frac{\Vert{\mathbf{E}}_{0}(S_{j})\Vert_{p}}{j^{3/2}}\Biggr)
^{p} .
\]
Theorem~\ref{directprop} then implies the following Burkholder-type inequality
that was established by Peligrad, Utev and Wu (\citeyear{PelUteWu07}),
Theorem 1:
\begin{corollary}
Let $p>2$ be a real number, and let $(X_{n})$ be an adapted stationary sequence
in the sense of Notation~\ref{notation1}. Then, for any integer $n$,
\[
{\mathbf{E}} \Bigl({\max_{1\leq j\leq n}}|S_{j}|^{p} \Bigr)\ll n^{p/2}{\mathbf{E}%
}(|X_{1}|^{p})+n^{p/2}\Biggl(\sum_{j=1}^{n}\frac{\Vert{\mathbf{E}}_{0}%
(S_{j})\Vert_{p}}{j^{3/2}}\Biggr)^{p} .
\]
\end{corollary}
\begin{pf*}{Proof of Lemma~\ref{comparison}}
We shall first prove that
for any positive integer~$k$,
%
%
\begin{equation}\label{comp1}%
\Vert{\mathbf{E}}_{0}(S_{2^{k}}^{2})\Vert_{p/2}\leq2^{k+1}\Vert{\mathbf
{E}%
}_{0}(X_{1}^{2})\Vert_{p/2}+2^{k+2}\Biggl( \sum_{j=0}^{k-1}\frac
{\Vert{\mathbf{E}}_{0}(S_{2^{j}})\Vert_{p}}{2^{j/2}}\Biggr) ^{2} .
\end{equation}
By using the notation $\bar{S}_{n}=X_{n+1}+\cdots+X_{n}$ and the fact that
$S_{2n}^{2}=S_{n}^{2}+\bar{S}_{n}^{2}+2S_{n}\bar{S}_{n}$, we get, by
stationarity, that
\[
\Vert{\mathbf{E}}_{0}(S_{2n}^{2})\Vert_{p/2}\leq2\Vert{\mathbf{E}}_{0}%
(S_{n}^{2})\Vert_{p/2}+2\Vert{\mathbf{E}}_{0} (S_{n}{\mathbf{E}}_{n}%
(\bar{S}_{n}))\Vert_{p/2} .
\]
Now, the Cauchy--Schwarz inequality, applied twice, gives
\[
\Vert{\mathbf{E}}_{0} (S_{n}{\mathbf{E}}_{n}(\bar{S}_{n}))\Vert
_{p/2}\leq\Vert{\mathbf{E}}_{0}^{1/2}(S_{n}^{2}){\mathbf{E}}_{0}%
^{1/2}({\mathbf{E}}_{n}^{2}(\bar{S}_{n}))\Vert_{p/2}\leq\Vert{\mathbf
{E}}%
_{0}(S_{n}^{2})\Vert_{p/2}^{1/2}\Vert{\mathbf{E}}_{0}(S_{n})\Vert_{p} .
\]
Hence, setting $b_{n}=\Vert{\mathbf{E}}_{0}(S_{n}^{2})\Vert_{p/2}$, it follows
that
\[
b_{2n}\leq2b_{n}+2b_{n}^{1/2}\Vert{\mathbf{E}}_{0}(S_{n})\Vert_{p} .
\]
By recurrence, this gives that
\[
b_{2^{k}}\leq2^{k}b_{0}+\sum_{j=0}^{k-1}2^{k-j}b_{2^{j}}^{1/2}\Vert
{\mathbf{E}}_{0}(S_{2^{j}})\Vert_{p} .
\]
With the notation $B_{k}=\max_{0\leq j\leq k}2^{-j}b_{2^{j}}$, we
derive that
\[
B_{k}\leq2\max\Biggl(b_{0},B_{k}^{1/2}\sum_{j=0}^{k-1}2^{-j/2}\Vert{\mathbf{E}
}_{0}(S_{2^{j}})\Vert_{p} \Biggr)
\]
implying that
\[
2^{-k}b_{2^{k}}\leq B_{k}\leq2b_{0}+2^{2} \Biggl(\sum_{j=0}^{k-1}2^{-j/2}%
\Vert{\mathbf{E}}_{0}(S_{2^{j}})\Vert_{p}\Biggr)^{2} .
\]
This ends the proof of inequality (\ref{comp1}).

We turn now to the proof of (\ref{comp2}). Let $r$ be the positive integer
such that $2^{r-1}\leq n<2^{r}$. Starting with the binary expansion
(\ref{binaryexpansion}), and using Minkowski's inequality twice, first with
respect to the conditional expectation, and second with respect to the
norm in
${\mathbf{L}}^{p}$, we get by stationarity that\looseness=-1
\[
\Vert{\mathbf{E}}_{0}(S_{n}^{2})\Vert_{p/2}\leq\Biggl(\sum_{k=0}^{r-1}%
b_{k}\Vert({\mathbf{E}}_{0}(S_{2^{k}}^{2}))^{1/2}\Vert_{p} \Biggr)^{2}%
\leq\Biggl(\sum_{k=0}^{r-1}\Vert{\mathbf{E}}_{0}(S_{2^{k}}^{2})\Vert
_{p/2}^{1/2} \Biggr)^{2} .
\]\looseness=0
Using then inequality (\ref{comp1}), we derive that
%
%
\begin{equation} \label{binary}%
\Vert{\mathbf{E}}_{0}(S_{n}^{2})\Vert_{p/2}\ll n\Vert{\mathbf{E}}_{0}%
(X_{1}^{2})\Vert_{p/2}+n\Biggl( \sum_{j=0}^{r-1}\frac{\Vert{\mathbf{E}}%
_{0}(S_{2^{j}})\Vert_{p}}{2^{j/2}}\Biggr) ^{2} .
\end{equation}
Since $(\Vert{\mathbf{E}}_{0}(S_{n})\Vert_{p})_{n\geq1}$ is
subadditive, using
item (1) of Lemma~\ref{lmasubadd}, inequality (\ref{comp2}) follows from
(\ref{binary}).
\end{pf*}

\subsection{Rosenthal inequalities for martingales and the case of even
powers}
\label{sectionMA}

\subsubsection{The martingale case}
\label{martingalecase}

For any real $p>2$, Theorem~\ref{directprop} applied to stationary martingale
differences gives the following inequality:
\[
{\mathbf{E}} \Bigl({\max_{1\leq j\leq n}}|S_{j}|^{p} \Bigr)\ll n{\mathbf{E}%
}(|X_{1}|^{p})+n \Biggl(\sum_{k=1}^{n}\frac{1}{k^{1+2\delta/p}}\Vert{\mathbf
{E}%
}_{0}(S_{k}^{2})\Vert_{p/2}^{\delta} \Biggr)^{p/(2\delta)} ,
\]
where $\delta=\min(1,1/(p-2))$.

Since for stationary martingale differences we have ${\mathbf{E}}(S_{n}%
^{2})=n{\mathbf{E}}(X_{1}^{2})$, we can express the inequality in the
following form useful for applications:
%
%
\begin{eqnarray}\label{consdirectmart}%
{\mathbf{E}} \Bigl({\max_{1\leq j\leq n}}|S_{j}|^{p} \Bigr)&\ll& n^{p/2}%
({\mathbf{E}}(X_{1}^{2}))^{p/2}+n{\mathbf{E}}(|X_{1}|^{p})\nonumber\\[-8pt]\\[-8pt]
&&{}+n \Biggl(\sum
_{k=1}^{n}\frac{1}{k^{1+2\delta/p}}\Vert{\mathbf{E}}_{0}(S_{k}^{2}%
)-{\mathbf{E}}(S_{k}^{2})\Vert_{p/2}^{\delta} \Biggr)^{p/(2\delta)}
.\nonumber
\end{eqnarray}
As we shall see in the next result, for a stationary sequence $(d_{i}%
)_{i\in{\mathbf{Z}}}$ of martingale differences in ${\mathbf{L}}^{p}$ with
$p\geq4$ an even integer,\vadjust{\goodbreak} this inequality can be sharpened since it
holds with
$\delta=2/(p-2)$ (see Comment~\ref{commentgreaterdelta}). As a
consequence, we
recover, in case $p=4$, inequality (1.6) stated in \citet{Rio09}
that was
obtained for variables in ${\mathbf{L}}^{q}$ with $q=p/2$, by using
Burkholder's classical inequality combined with Theorem 3 in
\citet{WuZha08}. Notice that inequality (1.6) stated in
\citet{Rio09} cannot be
generalized for $p>4$ since Theorem 3 in \citet{WuZha08} is only
valid for
variables in ${\mathbf{L}}^{q}$ with $1<q\leq2$.
\begin{theorem}
\label{mart2}Let $p\geq4$ be an even integer, and let $d_{0}$ be a real random
variable in ${\mathbf{L}}^{p}$, measurable with respect to $\mathcal{F}_{0}$
and such that ${\mathbf{E}}(d_{0}|{\mathcal{F}}_{-1})=0$. Let $d_{i}%
=d_{0}\circ T^{i}$ and $S_{n}=\sum_{i=1}^{n}d_{i}$. Then for any
integer $n$,
\[
{\mathbf{E}} \Bigl({\max_{1\leq j\leq n}}|S_{j}|^{p} \Bigr)\ll n {\mathbf{E}%
}(|d_{1}|^{p})+n \Biggl(\sum_{k=1}^{n}\frac{1}{k^{1+4/p(p-2)}}\Vert{\mathbf{E}
}_{0}(S_{k}^{2})\Vert_{p/2}^{2/(p-2)} \Biggr)^{p(p-2)/4} .
\]
\end{theorem}

The technique that makes this result possible is a special
symmetrization for
martingales initiated by \citet{KwaWoy91}.
\begin{proposition}
\label{tangent}Assume that $(e_{k})_{k}$ are stationary martingale differences
adapted to an increasing filtration $(\mathcal{F}_{k})_{k}$ that are
conditionally symmetric (the distribution of $e_{k}$ given $\mathcal{F}_{k-1}$
is equal to the distribution of $-e_{k}$ given $\mathcal{F}_{k-1}$). Assume,
in addition, that the $e_{k}$'s are conditionally independent given a sigma
algebra $\mathcal{G}$, and such that the law of $e_{k}$ given $\mathcal
{G}$ is
the same as the law of $e_{k}$ given $\mathcal{F}_{k-1}$. Let $S_{n}%
=\sum_{i=1}^{n}e_{i}$. Then for any even integer $p\geq4$ and any integer
$n\geq1$,
%
%
\begin{eqnarray} \label{ine2tangent}%
&&{\mathbf{E}} \Bigl({\max_{1\leq j\leq n}}|S_{j}|^{p} \Bigr)\nonumber\\[-8pt]\\[-8pt]
&&\qquad\ll n {\mathbf{E}%
}(|e_{1}|^{p})+n \Biggl(\sum_{k=1}^{n}\frac{1}{k^{1+4/p(p-2)}}\Vert{\mathbf{E}
}_{0}(S_{k}^{2})\Vert_{p/2}^{2/(p-2)} \Biggr)^{p(p-2)/4} .\nonumber
\end{eqnarray}
\end{proposition}
\begin{pf}
Due to Doob's\vspace*{1pt}
maximal inequality, we have that\break $\Vert{\max_{1\leq j\leq n}}|S_{j}|\Vert
_{p}\leq q\Vert S_{n}\Vert_{p}$ where $q=p(p-1)^{-1}$. Then it suffices to
show that inequality (\ref{ine2tangent}) holds for ${\mathbf{E}}%
(|S_{n}|^{p})$. We shall base this proof again on dyadic induction. Denote
$\bar{S}_{n}=$ $e_{n+1}+\cdots+e_{2n}$ and $a_{n}=\Vert S_{n}\Vert_{p}$.

We start from inequality (\ref{b1a2n}). Since the sequence of martingale
differences $(e_{k})$ is conditionally symmetric and conditionally independent
given a master sigma algebra $\mathcal{G}$, we have ${\mathbf{E}}(S_{n}%
^{p-1}{\bar{S}}_{n})+{\mathbf{E}}(S_{n}{\bar{S}}_{n}^{p-1})=0$ and therefore
\[
a_{2n}^{p}\leq2a_{n}^{p}+2^{p} \bigl({\mathbf{E}}(S_{n}^{p-2}\bar{S}_{n}%
^{2})+{\mathbf{E}}(\bar{S}_{n}^{p-2}S_{n}^{2}) \bigr) .
\]
Using Lemma~\ref{basic}, we have that
\[
{\mathbf{E}}(S_{n}^{p-2}\bar{S}_{n}^{2})\leq a_{n}^{p-2}\Vert{\mathbf
{E}}%
_{0}(S_{n}^{2})\Vert_{p/2}
\]
and
\[
{\mathbf{E}}(S_{n}^{2}\bar{S}_{n}%
^{p-2})\leq a_{n}^{p-4/(p-2)}\Vert{\mathbf{E}}_{0}(S_{n}^{2})\Vert
_{p/2}^{2/(p-2)} .
\]
Therefore, by combining all these bounds, we obtain for every even integer
$p\geq4$,
\begin{eqnarray*}
a_{2n}^{p} & \leq & 2a_{n}^{p}+2^{p} \bigl(a_{n}^{p-2}\Vert{\mathbf{E}}%
_{0}(S_{n}^{2})\Vert_{p/2}+a_{n}^{p-4/(p-2)}\Vert{\mathbf{E}}_{0}(S_{n}%
^{2})\Vert_{p/2}^{2/(p-2)} \bigr)\\
& \leq & 2a_{n}^{p}+2^{p+1}a_{n}^{p-4/(p-2)}\Vert{\mathbf{E}}_{0}(S_{n}%
^{2})\Vert_{p/2}^{2/(p-2)} .
\end{eqnarray*}
We end the proof by using Lemma~\ref{reclemma}.
\end{pf}
\begin{pf*}{Proof of Theorem~\ref{mart2}}
We consider our martingale differences sequence $(d_{k})_{k}$, and we construct
two decoupled tangent versions $(e_{k})_{k}$ and $(\tilde{e}_{k})_{k}$ to
$(d_{k})_{k}$ that are $\mathcal{G}$-conditionally independent between them
[here $\mathcal{G}=\sigma(\{d_{i}\})$]. This means that these two sequences
are martingale differences with the additional property that the conditional
distribution of $d_{k}$ given $\mathcal{F}_{k-1}$ is equal to the distribution
of $e_{k}$ given $\mathcal{F}_{k-1}$ and also to the distribution of
$\tilde{e}_{k}$ given $\mathcal{F}_{k-1}$; see Definition~6.1.4 in
\citet{delGin99} for the definition of a decoupled tangent
sequence, and their Proposition 6.1.5. for the crucial fact that decoupled
sequences always exist. We refer also to their Remark 6.1.6, for the
construction of $(e_{k})_{k}$ and~$(\tilde{e}_{k})_{k}$. Therefore, for any
even integer $p$,
\[
{\mathbf{E}} \Biggl(\sum_{i=1}^{n}d_{i} \Biggr)^{p}={\mathbf{E}} \Biggl(\sum
_{i=1}^{n}\bigl(d_{i}-{\mathbf{E}}(e_{i}|\mathcal{G}) \bigr)\Biggr)^{p}\leq{\mathbf{E}%
} \Biggl(\sum_{i=1}^{n}(d_{i}-e_{i}) \Biggr)^{p} .
\]
Now we use Corollary 6.6.8. in \citet{delGin99}; see also
\citet{Zin85}. Since $(e_{i}-\tilde{e}_{i})_{i}$ is a decoupled tangent
sequence of $(d_{i}-e_{i})_{i}$, it follows that
\[
{\mathbf{E}} \Biggl(\sum_{i=1}^{n}d_{i} \Biggr)^{p}\ll{\mathbf{E}} \Biggl(\sum
_{i=1}^{n}(e_{i}-\tilde{e}_{i}) \Biggr)^{p} .
\]
Notice that the distribution of $e_{i}-\tilde{e}_{i}$ is conditionally
symmetric given $\mathcal{G}$. Therefore, using Doob's maximal inequality
and applying Proposition~\ref{tangent}, we obtain that for every even integer
$p\geq4$ and any integer $n$,
%
%
\begin{eqnarray}\label{b1thmmart2}
&&{\mathbf{E}} \Biggl(\max_{1\leq k\leq n} \Biggl(\sum_{i=1}^{k}d_{i} \Biggr)^{p}%
\Biggr)\nonumber\\
&&\qquad\ll{\mathbf{E}} \Biggl(\sum_{i=1}^{n}(e_{i}-\tilde{e}_{i}) \Biggr)^{p}%
\nonumber\\[-8pt]\\[-8pt]
&&\qquad\ll n {\mathbf{E}}(|e_{1}-\tilde{e}_{1}|^{p})\nonumber\\
&&\qquad\quad{}+n \Biggl(\sum_{k=1}^{n}\frac
{1}{k^{1+4/p(p-2)}} \Biggl\|{\mathbf{E}}_{0} \Biggl( \Biggl(\sum_{i=1}^{k}%
(e_{i}-\tilde{e}_{i}) \Biggr)^{2} \Biggr) \Biggr\|_{p/2}^{2/(p-2)} \Biggr)^{p(p-2)/4}%
.\nonumber
\end{eqnarray}
Notice now that $\sum_{i=1}^{k}{\mathbf{E}}(e_{i}^{2}|{\mathcal{F}}%
_{i-1})=\sum_{i=1}^{k}{\mathbf{E}}(d_{i}^{2}|{\mathcal{F}}_{i-1})$
since both
quantities are obtained using only the conditional distributions of the
$d_{i}$'s and $e_{i}$'s, respectively, and these two sequences are tangent.
Tangency also implies that $d_{i}$ and $e_{i}$ have the same distributions.
Hence $\Vert d_{1}\Vert_{p}=\Vert e_{1}\Vert_{p}$. For the same
reasons, we
also have $\sum_{i=1}^{k}{\mathbf{E}}(\tilde{e}_{i}^{2}|{\mathcal{F}}%
_{i-1})=\sum_{i=1}^{k}{\mathbf{E}}(d_{i}^{2}|{\mathcal{F}}_{i-1})$ and
$\Vert
d_{1}\Vert_{p}=\Vert\tilde{e}_{1}\Vert_{p}$. Therefore, $\Vert
e_{1}-\tilde
{e}_{1}\Vert_{p}\leq2\Vert d_{1}\Vert_{p}$ and
\begin{eqnarray*}
\Biggl\|{\mathbf{E}}_{0} \Biggl( \Biggl(\sum_{i=1}^{k}(e_{i}-\tilde{e}%
_{i}) \Biggr)^{2} \Biggr) \Biggr\|_{p/2} & \leq & 2 \Biggl\|{\mathbf{E}}_{0}%
\Biggl(\sum_{i=1}^{k}{\mathbf{E}}(e_{i}^{2}|{\mathcal{F}}_{i-1}%
) \Biggr) \Biggr\|_{p/2}\\
&&{}+2 \Biggl\|{\mathbf{E}}_{0} \Biggl(\sum_{i=1}^{k}{\mathbf{E}%
}(\tilde{e}_{i}^{2}|{\mathcal{F}}_{i-1}) \Biggr) \Biggr\|_{p/2}\\
& = & 4 \Biggl\|{\mathbf{E}}_{0} \Biggl(\sum_{i=1}^{k}{\mathbf{E}}(d_{i}%
^{2}|{\mathcal{F}}_{i-1}) \Biggr) \Biggr\|_{p/2}\\
&=&4 \Biggl\|\sum_{i=1}^{k}%
{\mathbf{E}}_{0}(d_{i}^{2}) \Biggr\|_{p/2} .
\end{eqnarray*}
Theorem~\ref{mart2} follows by introducing these bounds in inequality
(\ref{b1thmmart2}).
\end{pf*}

\subsubsection{Application to stationary processes via martingale
approximation}
\label{sectionmartappro}

Theorem~\ref{mart2} together with the martingale
approximation provide an alternative Rosenthal-type inequality
involving the
projection operator, very useful for analyzing linear processes.

The next lemma is a slight reformulation of the martingale
approximation result
that can be found in the paper by Wu and Woodroofe (\citeyear{WuWoo04}),
Theorem~1; see also
\citet{ZhaWoo08} and \citet{GorPel11}.
\begin{lemma}
\label{mdec}Let $p\geq1$, and let $(X_{n})$ be an adapted stationary sequence
in the sense of Notation~\ref{notation1}. Then there is a triangular
array of
row-wise stationary martingale differences satisfying
%
%
\begin{equation}\label{defdiffmart}%
D_{0}^{n}=\frac{1}{n}\sum_{i=1}^{n}\bigl({\mathbf{E}}_{1}(S_{i})-{\mathbf{E}}
_{0}(S_{i})\bigr);\qquad D_{k}^{n}=D_{0}^{n}\circ T^{k},
\end{equation}
such that for any $1\leq k\leq n$ we have
%
%
\begin{equation}\label{MA}%
S_{k}=M_{k}^{n}+R_{k}^{n}\qquad\mbox{where }M_{k}^{n}=\sum_{i=1}^{k}D_{i}^{n}
\end{equation}
and
\[
{\max_{1\leq k\leq n}}\Vert R_{k}^{n}\Vert_{p}\leq2\Vert X_{0}\Vert
_{p}+\frac
{3}{n}\sum_{i=1}^{n}\Vert{\mathbf{E}}_{0}(S_{i})\Vert_{p} .
\]
\end{lemma}

We state now the Rosenthal-type inequality that we shall establish with the
help of the approximation result above.
\begin{theorem}
\label{general}Let $p\geq4$ be an even integer, and let $(X_{i})_{i\in
\mathbb{Z}}$ be as in Theorem~\ref{directprop}. Then the following inequality
is valid: for any integer $n$,
\begin{eqnarray*}
{\mathbf{E}} \Bigl({\max_{1\leq j\leq n}}|S_{j}|^{p} \Bigr) & \ll & n\Vert
X_{0}\Vert_{p}^{p}+n^{1-p} \Biggl(\sum_{i=1}^{n}\Vert{\mathbf{E}}_{0}%
(S_{i})\Vert_{p} \Biggr)^{p}\\
&&{} +n \Biggl(\sum_{k=1}^{n}\frac{1}{k^{1+4/p(p-2)}}\Vert{\mathbf{E}}_{0}%
(S_{k}^{2})\Vert_{p/2}^{2/(p-2)} \Biggr)^{p(p-2)/4} .
\end{eqnarray*}
\end{theorem}
\begin{remark}
\label{rmkcompmartge} Theorems~\ref{directprop} and~\ref{general} are, in
general, not comparable. Indeed, for $p\geq4$, Theorem~\ref{directprop} applies
with $\delta=1/(p-2)$, so the last term of the inequality stated in Theorem
\ref{general} can be bounded by the last term in the inequality from Theorem
\ref{directprop} [see item (2) of Comment~\ref{commentgreaterdelta}]. However,
the second term in Theorem~\ref{general} gives additional contribution.
\end{remark}
\begin{pf*}{Proof of Theorem~\ref{general}} We bound first $\max_{1\leq
k\leq n}{\mathbf{E}}(|S_{k}|^{p})$. By the martingale approximation of Lemma
\ref{mdec} combined with Theorem~\ref{mart2}, we get that
\begin{eqnarray*}
\max_{1\leq k\leq n}{\mathbf{E}}(|S_{k}|^{p})&\ll&\max_{1\leq k\leq
n}{\mathbf{E}}(|M_{k}^{n}|^{p})+\max_{1\leq k\leq n}{\mathbf{E}}(|R_{k}%
^{n}|^{p})\\
&\ll& n {\mathbf{E}}(|D_{0}^{n}|^{p})+\Vert X_{0}\Vert_{p}^{p}
+ \Biggl(\frac{1}{n}\sum_{i=1}^{n}\Vert{\mathbf{E}}_{0}(S_{i})\Vert
_{p} \Biggr)^{p}\\
&&{}+n \Biggl(\sum_{k=1}^{n}\frac{1}{k^{1+4/p(p-2)}}\Vert{\mathbf{E}%
}_{0}((M_{k}^{n})^{2})\Vert_{p/2}^{2/(p-2)} \Biggr)^{p(p-2)/4} .
\end{eqnarray*}
It remains to analyze the first and the last terms. By (\ref{MA})
applied with
$k=1$, we notice that $D_{1}^{n}=X_{1}-R_{1}^{n}$, and by the triangle
inequality,%
\[
\|D_{0}^{n}\|_{p}^{p}\leq(\|X_{0}\|_{p}+\|R_{1}^{n}\|_{p})^{p}\ll\|X_{0}
\|_{p}^{p}+\frac{1}{n^{p}} \Biggl(\sum_{i=1}^{n}\Vert{\mathbf{E}}_{0}%
(S_{i})\Vert_{p} \Biggr)^{p} .
\]
For analyzing the last term, we use the fact that
\[
{\mathbf{E}}_{0}((M_{k}^{n})^{2})\leq2{\mathbf{E}}_{0}(S_{k}^{2}%
)+2{\mathbf{E}}_{0}((R_{k}^{n})^{2}) .
\]
By using Lemma~\ref{mdec} it follows that
\[
\Vert{\mathbf{E}}_{0}((M_{k}^{n})^{2})\Vert_{p/2}\ll\Vert{\mathbf{E}}%
_{0}(S_{k}^{2})\Vert_{p/2}+\Vert X_{0}\Vert_{p}^{2}+n^{-2} \Biggl(\sum_{i=1}%
^{n}\Vert{\mathbf{E}}_{0}(S_{i})\Vert_{p} \Biggr)^{2}
\]
and overall
\begin{eqnarray*}
\max_{1\leq k\leq n}{\mathbf{E}}(|S_{k}|^{p}) & \leq & n\Vert X_{0}\Vert
_{p}^{p}+n^{1-p} \Biggl(\sum_{i=1}^{n}\Vert{\mathbf{E}}_{0}(S_{i})\Vert
_{p} \Biggr)^{p}\\
&&{} + n \Biggl(\sum_{k=1}^{n}\frac{1}{k^{1+4/p(p-2)}}\Vert{\mathbf{E}}_{0}%
(S_{k}^{2})\Vert_{p/2}^{2/(p-2)} \Biggr)^{p(p-2)/4} .
\end{eqnarray*}
We apply\vspace*{1pt} now inequality (\ref{cons1sub}) that has the effect of an addition
of a forth term, namely $n (\sum_{k=1}^{n}k^{-1-1/p}%
\Vert{\mathbf{E}}_{0}(S_{k})\Vert_{p} )^{p}$. However, we can express our
inequality without including this term because it can be bounded, up to a
multiplicative constant, by $n^{1-p} (\sum_{i=1}^{n}\Vert{\mathbf{E}}%
_{0}(S_{i})\Vert_{p} )^{p}$. Indeed, notice that it is enough to show that
for a certain universal constant $C$,
\[
\max_{1\leq i\leq n}\Vert{\mathbf{E}}_{0}(S_{i})\Vert_{p}\leq\frac{C}{n}
\sum_{i=1}^{n}\Vert{\mathbf{E}}_{0}(S_{i})\Vert_{p} .
\]
To prove it, we first notice that
\[
{\max_{1\leq i\leq n}}\Vert{\mathbf{E}}_{0}(S_{i})\Vert_{p}\leq
{\max_{1\leq
i\leq[ n/2]}}\Vert{\mathbf{E}}_{0}(S_{i})\Vert_{p}+{\max_{[n/2]<i\leq
n}}\Vert{\mathbf{E}}_{0}(S_{i})\Vert_{p}%
\]
and that for any $i\in\{1,\ldots,[n/2]\}$,
\[
\Vert{\mathbf{E}}_{0}(S_{i})\Vert_{p}\leq\bigl\Vert{\mathbf
{E}}_{0}\bigl(S_{i+[n/2]}%
\bigr)\bigr\Vert_{p}+\bigl\Vert{\mathbf{E}}_{0}\bigl(S_{i+[n/2]}-S_{i}\bigr)
\bigr\Vert_{p} .
\]
Therefore, by the properties of conditional expectation and
stationarity, it
follows that%
\[
\max_{1\leq i\leq n}\Vert{\mathbf{E}}_{0}(S_{i})\Vert_{p}\leq\bigl\Vert
{\mathbf{E}%
}_{0}\bigl(S_{[n/2]}\bigr)\bigr\Vert_{p}
+2\max_{[n/2]<i\leq n}\Vert{\mathbf{E}}_{0}%
(S_{i})\Vert_{p} .
\]
To complete the proof, it remains to apply inequality (\ref{condsubadd2})
to the subadditive sequence $(\Vert{\mathbf{E}}_{0}(S_{i})\Vert
_{p})_{i\geq1}%
$.
\end{pf*}

\subsection{Rosenthal inequality in terms of individual summands}
\label{sectionindsum}

For the sake of applications in this section we indicate
how to estimate the terms that appear in our Rosenthal inequalities in terms
of individual summands and formulate some specific inequalities. By
substracting ${\mathbf{E}}(S_{k}^{2})$ and applying the triangle inequality
we can reformulate all the inequalities in terms of the quantities
${\mathbf{E}}(S_{k}^{2})$, $\|{\mathbf{E}}_{0}(S_{k})\|_{p}$ and
$\Vert{\mathbf{E}}_{0}(S_{k}^{2})-{\mathbf{E}}(S_{k}^{2})\Vert_{p/2}$.
The next
lemma proposes a simple way to estimate these quantities in terms of
coefficients in the spirit of \citet{Gor69}.
\begin{lemma}
\label{mixing}Under the stationary setting assumptions in Notation
\ref{notation1}, we have the following estimates:
%
%
\begin{eqnarray}\label{var}%
{\mathbf{E}}(S_{k}^{2})&\leq&2k\sum_{j=0}^{k-1}|{\mathbf{E}}(X_{0}X_{j})| ,
\\
%
%
\label{condterm}%
\|{\mathbf{E}}_{0}(S_{k})\|_{p}&\leq&\sum_{\ell=1}^{n}\|{\mathbf{E}}_{0}%
(X_{\ell})\|_{p}
\end{eqnarray}
and
%
%
\begin{eqnarray}\label{varsquare}
&&\Vert{\mathbf{E}}_{0}(S_{k}^{2})-{\mathbf{E}}(S_{k}^{2})\Vert_{p/2}\nonumber\\
&&\qquad\leq
2\sum_{i=1}^{k}\sum_{j=0}^{k-i}\Vert{\mathbf{E}}_{0}(X_{i}X_{i+j}%
)-{\mathbf{E}}(X_{i}X_{i+j})\Vert_{p/2}\nonumber\\[-8pt]\\[-8pt]
&&\qquad\leq2\sum_{i=1}^{k}\sum_{j=0}^{k-i}\sup_{\ell\geq0}\Vert{\mathbf{E}}_{0}
(X_{i}X_{i+\ell})-{\mathbf{E}}(X_{i}X_{i+\ell})\Vert_{p/2}\wedge(2\Vert
X_{0}{\mathbf{E}}_{0}(X_{j})\Vert_{p/2})\nonumber\\
&&\qquad\leq4\sum_{j=1}^{k}j\Vert X_{0}{\mathbf{E}}_{0}(X_{j})\Vert_{p/2}+2\sum
_{i=1}^{k}i\sup_{j\geq i}\Vert{\mathbf{E}}_{0}(X_{i}X_{j})-{\mathbf{E}}%
(X_{i}X_{j})\Vert_{p/2}.\nonumber
\end{eqnarray}
\end{lemma}

Mixing coefficients are useful to continue the estimates from Lemma
\ref{mixing}. We refer to the books by Bradley [(\citeyear{Bra07}),
Theorem 4.13 via Remark 4.7, item VI], Rio [(\citeyear{Rio00}),
Theorem 2.5 and Appendix, Section C] and Dedecker
et al. [(\citeyear{Dedetal07}), Remark 2.5 and Chapter 3] for various
estimates of the
coefficients involved in Lemma~\ref{mixing} and examples. We shall also
provide applications and explicit computations of the quantities involved.

We formulate the following proposition:
\begin{proposition}
\label{consdirect} Let $p>2$ be a real number, and let $(X_{i})_{i\in
\mathbf{Z}}$ be a stationary sequence of real-valued random variables
in ${\mathbf{L}}_{p}$ adapted to an increasing
filtration~$({\mathcal{F}}_{i})$. For any $j\geq1$, let
%
%
\begin{equation}\label{notalambda}%
\lambda(j)=\max\Bigl(\Vert X_{0}{\mathbf{E}}_{0}(X_{j})\Vert_{p/2},\sup
_{i\geq
j}\Vert{\mathbf{E}}_{0}(X_{i}X_{j})-{\mathbf{E}}(X_{i}X_{j})\Vert
_{p/2} \Bigr) .
\end{equation}
Then for every positive integer $n$,
\begin{eqnarray*}
\Bigl\|{\max_{1\leq j\leq n}}|S_{j}|\Bigr\|_{p}&\ll& n^{1/2}%
\Biggl(\sum_{k=0}^{n-1}|{\mathbf{E}}(X_{0}X_{k})| \Biggr)^{1/2}\\
&&{}+n^{1/p}\|X_{1}\|_{p}
+cn^{1/p}\sum_{k=1}^{n}\frac{1}{k^{1/p}}\|{\mathbf{E}}_{0}(X_{k}%
)\|_{p}\\
&&{}+n^{1/p} \Biggl(\sum_{k=1}^{n}\frac{1}{k^{(2/p)-1}}(\log k)^{\gamma
}\lambda(k) \Biggr)^{1/2} ,
\end{eqnarray*}
where $\gamma$ can be taken $\gamma=0$ for $2<p\leq3$ and $\gamma>p-3$ for
$p>3;$ $c=1$ for $2<p<4$ and $c=0$ for $p\geq4$. The constant that is
implicitly involved in the notation $\ll$ depends on $p$ and $\gamma$,
but it
does not depend on $n$ and on the~$X_{i}$'s.
\end{proposition}
\begin{pf}
The proof of this
proposition is basically a combination of Theorem~\ref{directprop} and Lemma
\ref{mixing}. By the triangle inequality,
\[
\|{\mathbf{E}}_{0}(S_{k}^{2})\|_{p/2}\leq\Vert{\mathbf{E}}_{0}(S_{k}%
^{2})-{\mathbf{E}}(S_{k}^{2})\Vert_{p/2}+{\mathbf{E}}_{0}(S_{k}^{2}) .
\]
By (\ref{var}), for any $p>2$ and any $\delta>0$, we easily obtain
\[
\Biggl( \sum_{k=1}^{n}\frac{1}{k^{1+2\delta/p}}({\mathbf{E}}(S_{k}%
^{2}))^{\delta}\Biggr) ^{1/(2\delta)}\ll n^{1/2-1/p} \Biggl(\sum_{j=0}%
^{n-1}|{\mathbf{E}}(X_{0}X_{j})| \Biggr)^{1/2} .
\]
Then we use inequality (\ref{condterm}) and changing the order of summation
\[
\sum_{k=1}^{n}\frac{1}{k^{1+1/p}}\|{\mathbf{E}}_{0}(S_{k})\|_{p}\ll\sum
_{k=1}^{n}\frac{1}{k^{1/p}}\|{\mathbf{E}}_{0}(X_{k})\|_{p} .
\]
Now for the situation $0<\delta<1$, by H\"{o}lder's inequality,
\begin{eqnarray*}
&&\Biggl( \sum_{k=1}^{n}\frac{1}{k^{1+2\delta/p}}\Vert{\mathbf{E}}_{0}(S_{k}%
^{2})-{\mathbf{E}}(S_{k}^{2})\Vert_{p/2}^{\delta}\Biggr) ^{p/(2\delta)}%
\\
&&\qquad\ll\Biggl( \sum_{k=1}^{n}\frac{(\log k)^{\gamma}}{k^{1+2/p}}\Vert
{\mathbf{E}}_{0}(S_{k}^{2})-{\mathbf{E}}(S_{k}^{2})\Vert_{p/2}\Biggr)
^{p/2} ,
\end{eqnarray*}
where $\gamma>1/\delta-1$. We continue the proof by using (\ref{varsquare})
and get
\[
\sum_{k=1}^{n}\frac{(\log k)^{\gamma}}{k^{1+2/p}}\Vert{\mathbf{E}}_{0}%
(S_{k}^{2})-{\mathbf{E}}(S_{k}^{2})\Vert_{p/2}\ll\sum_{k=1}^{n}\frac
{(\log
k)^{\gamma}}{k^{(2/p)-1}}\lambda(k) .
\]
Proposition~\ref{consdirect} follows by using Theorem~\ref{directprop}
combined with all the above upper bounds.\vadjust{\goodbreak}
\end{pf}

We give now a consequence of Theorem~\ref{directprop} that will be used
in one
of our applications. The proof is omitted since it is in the spirit of the
proof of Proposition~\ref{consdirect}; namely, we use Lemma~\ref{mixing}
combined with H\"{o}lder's inequality and the fact that $\|{\mathbf{E}}%
_{0}(S_{k})\|_{p}\leq\|{\mathbf{E}}_{0}(S_{k}^{2})-{\mathbf{E}}(S_{k}%
^{2})\|_{p/2}^{1/2}+ ({\mathbf{E}}(S_{k}^{2}))^{1/2}$.
\begin{proposition}
\label{consdirect2} Let $p>2$ be a real number, and let $(X_{i})_{i\in
\mathbf{Z}}$ be a stationary sequence of real-valued random variables in
${\mathbf{L}}_{p}$ adapted to an increasing filtration~$({\mathcal{F}}_{i})$.
Let $(\lambda(j))_{j \geq1}$ be defined by (\ref{notalambda}). For every
positive integer $n$, the following inequality holds: for any
$\varepsilon
>0$,
\begin{eqnarray*}
{\mathbf{E}} \Bigl({\max_{1\leq j\leq n}}|S_{j}|^{p} \Bigr)&\ll& n^{p/2}%
\Biggl(\sum_{k=0}^{n-1}|{\mathbf{E}}(X_{0}X_{k})| \Biggr)^{p/2}%
+n{\mathbf{E}} ( |X_{1}|^{p})\\
&&{}+n \sum_{k=1}^{n}k^{p-2+\varepsilon}\lambda
^{p/2}(k).
\end{eqnarray*}
The constant that is implicitly involved in the notation $\ll$ depends
on $p$ and~$\varepsilon$, but it does not depend on $n$.
\end{proposition}

\section{Applications and examples}
\label{sectionappliexamples}

As we have seen, Propositions~\ref{maxinequality} and~\ref{propmaxineproba}
give a direct approach to compare the moments of order $p$ of the
maximum of
the partial sums to the moments of order $p$ of the partial sum. We
start this
section by presenting two additional applications of these propositions
to the
convergence of maximum of partial sums and to the maximal Bernstein inequality
for dependent structures. In the last three examples, we apply our
results on
the Rosenthal-type inequalities to different classes of processes. In
all the
examples given in Sections~\ref{sectionarch},~\ref{sectionFLP} and
\ref{sectionRMC}, we show that, for real numbers $p$ larger than $2$, the
order of magnitude of $\Vert{\max_{1\leq k\leq n}}|S_{k}|\Vert_{p}$ is
essentially given by the order of magnitude of $\Vert S_{n}\Vert_{2}$
(or by a
bound of it).

\subsection{Convergence of the maximum of partial sums in ${\mathbf{L}}^{p}$}

\begin{corollary}
\label{cortightlp} Let $p>1$, and let $(X_{i})_{i\in{\mathbf{Z}}}$ be a
strictly stationary sequence of centered real-valued random variables in
${\mathbf{L}}^{p}$ adapted to an increasing and stationary filtration
$({\mathcal{F}}_{i})_{i\in{\mathbf{Z}}}$. Assume that
%
%
\begin{equation} \label{hypoappmart}%
\lim_{n\rightarrow\infty}n^{-1/p}\Vert S_{n}\Vert_{p}=0 .
\end{equation}
Assume, in addition, that
%
%
\begin{equation} \label{maxwoodLp}%
\sum_{n\geq1}\frac{\Vert{\mathbf{E}}(S_{n}|{\mathcal{F}}_{0})\Vert_{p}%
}{n^{1+1/p}}<\infty.
\end{equation}
Then
%
%
\begin{equation} \label{resmaxtightlp}%
\lim_{n\rightarrow\infty}n^{-1/p}\Bigl\Vert{\max_{1\leq k\leq n}}
|S_{k}|\Bigr\Vert_{p}=0 .
\end{equation}
\end{corollary}
\begin{remark}
This corollary is particularly useful for studying the asymptotic
behavior of
a partial sum via a martingale approximation. Assume there exists a strictly
stationary sequence $(d_{i})_{i\in{\mathbf{Z}}}$ of martingale differences
with respect to $({\mathcal{F}}_{i})_{i\in{\mathbf{Z}}}$ that are in
${\mathbf{L}}^{p}$, such that $\lim_{n\rightarrow\infty}n^{-1/p}\Vert
S_{n}-\sum_{i=1}^{n}d_{i}\Vert_{p}=0$. Then, if the condition
(\ref{maxwoodLp}) holds for the sequence $(X_{i})_{i\in{\mathbf{Z}}}$,
by a
construction in \citet{ZhaWoo08} and by the uniqueness of the
martingale approximation, the sequence $(X_{i}-d_{i})_{i\in{\mathbf
{Z}}}$ is
still a strictly stationary sequence and by our theorem $\lim
_{n\rightarrow
\infty}n^{-1/p}\Vert{\max_{1\leq k\leq n}}|S_{k}-\sum_{i=1}^{k}d_{i}|\Vert
_{p}=0$. As a matter of fact, for $p=2$, our corollary leads to the functional
form of the central limit theorem for $\{n^{-1/2}S_{[nt]},t\in[0,1]\}$;
see also Theorem 1.1 in \citet{PelUte05}.
\end{remark}
\begin{pf*}{Proof of Corollary~\ref{cortightlp}}
Let $m$ be an integer and $k=k_{n,m}=[n/m]$ (where $[x]$ denotes the integer
part of $x$).

The initial step of the proof is to divide the variables in blocks of
size $m$
and to make the sums in each block. Let
\[
X_{i,m}=\sum_{j=(i-1)m+1}^{im}X_{j},\qquad i\geq1.
\]
Notice first that
\begin{eqnarray*}
\Biggl\Vert\sup_{t\in[0,1]} \Biggl|\sum_{j=1}^{[nt]}X_{j}-\sum_{i=1}%
^{[kt]}X_{i,m} \Biggr|\Biggr\Vert_{p}&\leq&\Biggl\Vert\sup_{t\in[0,1]} \Biggl|\sum
_{i=[kt]m+1}^{[nt]}X_{i} \Biggr|\Biggr\Vert_{p}\\
&\leq& m\Bigl\Vert\max_{1\leq i\leq n}%
X_{i}\Bigr\Vert_{p} .
\end{eqnarray*}
Since for every $\varepsilon>0$,
\[
{\mathbf{E}}\Bigl({\max_{1\leq i\leq n}}|X_{i}|^{p}\Bigr)\leq\varepsilon^{p}+\sum
_{i=1}^{n}{\mathbf{E}} \bigl(|X_{i}|^{p}\mathbf{1}_{\{|X_{i}|>\varepsilon
\}} \bigr) ,
\]
and since $\Vert X_{i}\Vert_{p}<\infty$ for all $i$, we derive that
$\lim_{n\rightarrow\infty}\Vert{\max_{1\leq i\leq n}X_{i}}\Vert_{p}/n^{1/p}=0$.
Hence, in order to prove (\ref{resmaxtightlp}), it remains to show that
%
%
\begin{equation} \label{butresmaxtightlp}%
\lim_{m\rightarrow\infty}\limsup_{n\rightarrow\infty}%
n^{-1/p}\Biggl\Vert\sup_{t\in[0,1]} \Biggl|\sum_{i=1}^{[kt]}X_{i,m}%
\Biggr|\Biggr\Vert_{p}=0 .
\end{equation}
Applying Proposition~\ref{maxinequality} to the variables
$(X_{i,m})_{1\leq
i\leq k}$ which are adapted with respect to ${\mathcal{F}}_{im}$, and taking
into account Remark~\ref{Rem1inegalite}, we get that
\[
\Biggl\Vert\sup_{t\in[0,1]} \Biggl|\sum_{i=1}^{[kt]}X_{i,m} \Biggr|\Biggr\Vert_{p}%
\ll\max_{1\leq j\leq k}\Biggl\Vert\sum_{\ell=1}^{jm}X_{\ell}\Biggr\Vert_{p}+k^{1/p}%
\sum_{j=1}^{k}\frac{\|{\mathbf{E}}(S_{jm}|{\mathcal
{F}}_{0})\|_{p}}{j^{1+1/p}%
} ,
\]
where for the last term we used the fact that for any positive integer $u$,
$\Vert{\mathbf{E}} (\sum_{j=um2^{\ell}+1}^{(u+1)m2^{\ell}}X_{j}%
|{\mathcal{F}}_{um2^{\ell}} )\Vert_{p}=\Vert{\mathbf{E}} (\sum
_{j=1}^{m2^{\ell}}X_{j}|{\mathcal{F}}_{0} )\Vert_{p}$. Condition
(\ref{hypoappmart}) implies that
\[
\max_{1\leq j\leq k}\Biggl\Vert\sum_{\ell=1}^{jm}X_{\ell}\Biggr\Vert
_{p}=o((km)^{1/p}%
)=o(n^{1/p}) .
\]
Now, by the subadditivity of the sequence $ (\|{\mathbf{E}}(S_{n}%
|{\mathcal{F}}_{0})\|_{p})_{n\geq1}$ and applying Lem\-ma~\ref{lmasubadd}
we have
%
%
\begin{equation} \label{fatou}\qquad
n^{-1/p}k^{1/p}\sum_{j=1}^{k}\frac{\|{\mathbf{E}}(S_{jm}|{\mathcal{F}}%
_{0})\|_{p}}{j^{1+1/p}}\ll\sum_{\ell=1}^{m}\frac{\|{\mathbf{E}}(S_{\ell
}|{\mathcal{F}}_{0})\|_{p}}{(\ell+m)^{1+1/p}}+\sum_{j\geq m}\frac
{\|{\mathbf{E}}(S_{j}|{\mathcal{F}}_{0})\|_{p}}{j^{1+1/p}} .
\end{equation}
Hence, by (\ref{maxwoodLp}) and by using the dominated convergence
theorem for
discrete measures applied to the first term in the right-hand side of
(\ref{fatou}), we get that
\[
\lim_{m\rightarrow\infty}\limsup_{n\rightarrow\infty}%
n^{-1/p}k^{1/p}\sum_{j=1}^{k}\frac{\|{\mathbf{E}}(S_{jm}|{\mathcal{F}}%
_{0})\|_{p}}{j^{1+1/p}}=0 ,
\]
which ends the proof of the corollary.
\end{pf*}

\subsection{Maximal exponential inequalities for strong mixing
sequences}

Let us first recall the definition of strongly mixing sequences,
introduced by
\citet{Ros56}: for any two $\sigma$ algebras $\mathcal{A}$ and
$\mathcal{B}$, we define the $\alpha$-mixing coefficient by
\[
\alpha(\mathcal{A},\mathcal{B})={\sup_{A\in\mathcal{A},B\in\mathcal{B}%
}}|\mathbf{P}(A\cap B)-\mathbf{P}(A)\mathbf{P}(B)|.
\]

Let $(X_{k},k\geq1)$ be a sequence of real-valued random variables
defined on
$( \Omega,\mathcal{A},\mathbf{P}) $. This sequence will be called
strongly mixing if
%
%
\begin{equation} \label{defalpha}%
\alpha(n):=\sup_{k\geq1}\alpha( \mathcal{F}_{k},\mathcal{G}%
_{k+n}) \rightarrow0\qquad\mbox{as }n\rightarrow\infty,
\end{equation}
where $\mathcal{F}_{j}:=\sigma(X_{i},i\leq j)$ and $\mathcal{G}_{j}%
:=\sigma(X_{i},i\geq j)$ for $j\geq1$.

In 2009, Merlev\`{e}de, Peligrad and Rio proved (see their Theorem 2)
that for a strongly mixing sequence of centered random variables
satisfying\break
${\sup_{i\geq1}}\Vert X_{i}\Vert_{\infty}\leq M$ and for a certain $c>0$
%
%
\begin{equation} \label{alphacond}%
\alpha(n)\leq\exp(-cn) ,
\end{equation}
the following Bernstein-type inequality is valid: there is a constant $C$
depending only on $c$ such that for all $n\geq2$,
%
%
\begin{equation}\label{bernSn}%
\mathbf{P}(|S_{n}|\geq x)\leq\exp\biggl(-{\frac{Cx^{2}}{v^{2}n+M^{2}+xM(\log
n)^{2}}} \biggr) ,
\end{equation}
where
%
%
\begin{equation} \label{defv2}%
v^{2}=\sup_{i>0} \biggl(\operatorname{Var}(X_{i})+2\sum_{j>i}|{\operatorname
{Cov}}%
(X_{i},X_{j})|\biggr) .
\end{equation}
Proving the maximal version of inequality (\ref{bernSn}) cannot be handled
directly using either Theorem 2.2 in \citet{MorSerSto82}
or using Theorem 1 in \citet{KevMas11} since the right-hand side of
(\ref{bernSn}) does not satisfy the assumptions of these papers.
However, an
application of Proposition~\ref{propmaxineproba} leads to the maximal version
of Theorem 2 in \citet{MerPelRio09}.
\begin{corollary}
\label{cormaxbernstein} Let $(X_{j})_{j\geq1}$ be a sequence of centered
real-valued \mbox{random} variables. Suppose that there exists a positive $M$ such
that ${\sup_{i\geq1}}\Vert X_{i}\Vert_{\infty}\leq M$ and that the strongly
mixing coefficients $(\alpha(n))_{n\geq1}$ of the sequence satisfy
(\ref{alphacond}). Then there exist constants $C=C(c)$ and $K=K(M,c)$ such
that for all integers $n\geq2$ and all real $x>K\log n$,
%
%
\begin{equation} \label{bern2}%
\mathbf{P}\Bigl({\max_{1\leq k\leq n}}|S_{k}|\geq x\Bigr)\leq\exp\biggl(-{\frac{Cx^{2}%
}{v^{2}n+M^{2}+xM(\log n)^{2}}} \biggr) .
\end{equation}
\end{corollary}
\begin{pf}
We first apply
inequality (\ref{inegalitemax2proba}) of Proposition~\ref{propmaxineproba}
with $p=2$ and $\varphi(x)=e^{t|x|}$ where $t$ is a positive integer.
According to Theorem 2 in \citet{MerPelRio09}, there exist
positive constants $C_{1}$ and $C_{2}$ depending only on $c$ such that
for all
$n\geq2$ and any positive $t$ such that $t<\frac{1}{C_{1}M(\log
n)^{2}}$, the
following inequality holds:
\[
\log\mathbf{E} (\exp(tS_{n}) )\leq\frac{C_{2}t^{2}(nv^{2}+M^{2}%
)}{1-C_{1}tM(\log n)^{2}} .
\]
Then an optimization on $t$ gives that there is a constant $C_{3}$ depending
only on $c$ such that for all $n\geq2$ and any positive real $x$,
%
%
\begin{eqnarray}\label{bern2*}%
&&\mathbf{P}\Bigl({\max_{1\leq k\leq n}}|S_{k}|\geq4x\Bigr)\nonumber\\
&&\qquad\leq\exp\biggl(-{\frac
{C_{3}x^{2}%
}{v^{2}n+M^{2}+xM(\log n)^{2}}} \biggr)\\
&&\qquad\quad{}+4x^{-2}\Biggl(\sum_{l=0}^{r-1} \Biggl(\sum_{k=1}^{2^{r-l}-1}\bigl\|{\mathbf{E}%
}\bigl(S_{v+(k+1)2^{l}}-S_{v+k2^{l}}|{\mathcal{F}}_{k2^{l}}\bigr)\bigr\|_{2}^{2}%
\Biggr)^{1/2} \Biggr)^{2} ,\nonumber
\end{eqnarray}
where $r$ is the positive integer satisfying $2^{r-1}<n\leq2^{r}$ and
$v=[x/M]$.

It remains to bound up the second term in the right-hand side of the above
inequality. Notice first that for any centered variable $Z$ such that
$\Vert
Z\Vert_{\infty}\leq B$, Ibragimov's covariance inequality [see Theorem
1.11 in
\citet{Bra07}] gives
\[
\Vert{\mathbf{E}}(Z|{\mathcal{F}})\Vert_{2}^{2}= \operatorname
{Cov}({\mathbf
{E}%
}(Z|{\mathcal{F}}),Z) \leq4B^{2}\alpha({\mathcal{F}},\sigma(Z)) .
\]
Therefore, applying this last estimate with $Z=S_{v+(k+1)2^{l}}-S_{v+k2^{l}}$
and ${\mathcal{F}}={\mathcal{F}}_{k2^{l}}$, we get that
\[
\bigl\|{\mathbf{E}}\bigl(S_{v+(k+1)2^{l}}-S_{v+k2^{l}}|{\mathcal
{F}}_{k2^{l}}\bigr)\bigr\|_{2}%
^{2}\leq4M^{2}2^{2l}\alpha(v)
\]
implying that
\[
\Biggl(\sum_{l=0}^{r-1} \Biggl(\sum_{k=1}^{2^{r-l}-1}\bigl\|{\mathbf{E}}%
\bigl(S_{v+(k+1)2^{l}}-S_{v+k2^{l}}|{\mathcal{F}}_{k2^{l}}\bigr)\bigr\|_{2}^{2}%
\Biggr)^{1/2} \Biggr)^{2}\leq4M^{2}2^{2r}\bigl(\sqrt{2}+1\bigr)^{2}\alpha(v) .
\]
Since $2^{2r}\leq4n^{2}$ and $[x/M]\geq x/(2M)$, for $x\geq2M$, by using
(\ref{alphacond}), we get that, for any $x\geq2M$,
%
%
\begin{eqnarray} \label{bern2*b2}%
&&\Biggl(\sum_{l=0}^{r-1} \Biggl(\sum_{k=1}^{2^{r-l}-1}\bigl\|{\mathbf{E}}%
\bigl(S_{v+(k+1)2^{l}}-S_{v+k2^{l}}|{\mathcal{F}}_{k2^{l}}\bigr)\bigr\|_{2}^{2}%
\Biggr)^{1/2} \Biggr)^{2}\nonumber\\[-8pt]\\[-8pt]
&&\qquad\leq3(2^{5})M^{2}n^{2}\exp\bigl(-cx/(2M)\bigr) .\nonumber
\end{eqnarray}
Starting from (\ref{bern2*}) and using (\ref{bern2*b2}), we then derive that
for any $x\geq2M\max(1$, $4c^{-1}\log n)$,
\[
{\mathbf{P}}(S_{2^{r}}^{\ast}\geq4x)\leq\exp\biggl(-{\frac{C_{3}x^{2}}%
{v^{2}n+M^{2}+xM(\log n)^{2}}} \biggr)+96\exp\biggl(-\frac{xc}{4M} \biggr)
\]
proving inequality (\ref{bern2}).
\end{pf}

\subsection{Application to ARCH models}
\label{sectionarch}

Theorem~\ref{directprop} applies to the case where $(X_{i})_{i\in
{\mathbf{Z}}%
}$ has an ARCH($\infty$) structure as described by \citet
{GirKokLei00}, that is,
%
%
\begin{equation}\label{defARCH}%
X_{n}=\sigma_{n}\eta_{n}\qquad \mbox{with } \sigma_{n}^{2}=c+\sum
_{j=1}^{\infty
}c_{j}X_{n-j}^{2} ,
\end{equation}
where $(\eta_{n})_{n\in\mathbf{Z}}$ is a sequence of i.i.d. centered random
variables such that \mbox{${\mathbf{E}}(\eta_{0}^{2})=1$}, and where $c\geq0$,
$c_{j}\geq0$ and $\sum_{j\geq1}c_{j}<1$. Notice that $(X_{i})_{i
\in{\mathbf{Z}}}$ is a stationary sequence of martingale differences adapted
to the filtration $({\mathcal{F}}_{i})$ where ${\mathcal{F}}_{i}= \sigma(
\eta_{k} , k \leq i)$.

Let $p>2$ and assume that $\Vert\eta_{0}\Vert_{p}<\infty$. Notice first that
%
%
\begin{equation}\label{ARCHfirstequality}%
\bigl\Vert{\mathbf{E}}(X_{j}^{2}|{\mathcal{F}}_{0})-{\mathbf{E}}(X_{0}^{2}%
)\bigr\Vert_{p/2}=\bigl\Vert{\mathbf{E}}(\sigma_{j}^{2}|{\mathcal
{F}}_{0})-{\mathbf{E}%
}(\sigma_{j}^{2})\bigr\Vert_{p/2} .
\end{equation}
In addition, since ${\mathbb{E}}(\eta_{0}^{2})=1$ and $\sum_{j\geq1}c_{j}<1$,
the unique stationary solution to (\ref{defARCH}) is given by
\citet{GirKokLei00},
%
%
\begin{equation} \label{solARCH}%
\sigma_{n}^{2}=c+c\sum_{\ell=1}^{\infty}\sum_{j_{1},\ldots,j_{\ell
}=1}^{\infty
}c_{j_{1}}\cdots c_{j_{\ell}}\eta_{n-j_{1}}^{2}\cdots\eta_{n-(j_{1}%
+\cdots+j_{\ell})}^{2} .
\end{equation}
Starting from (\ref{ARCHfirstequality}) and using (\ref{solARCH}), one can
prove that
\[
\bigl\Vert{\mathbf{E}}(X_{j}^{2}|{\mathcal{F}}_{0})-{\mathbf{E}}(X_{0}^{2}%
)\bigr\Vert_{p/2}\leq2c\Vert\eta_{0}\Vert_{p}^{2}\sum_{\ell=1}^{\infty}\ell
\kappa^{\ell-1}\sum_{i=[j/\ell]}^{\infty}c_{i} ,
\]
where $\kappa=\Vert\eta_{0}\Vert_{p}^{2}\sum_{j\geq1}c_{j}$; see
Section 6.6
in \citet{DedMer11} for more detailed computations.
Therefore, if
%
%
\begin{equation} \label{condARCH}%
\Vert\eta_{0}\Vert_{p}^{2}\sum_{j\geq1}c_{j}<1
\end{equation}
for any $p>2$ and $\delta\in\ ]0,1]$, we obtain%
\[
\sum_{k=1}^{n}\frac{1}{k^{1+2\delta/p}}\|{\mathbf{E}}_{0}(S_{k}^{2}%
)-{\mathbf{E}}(S_{k}^{2})\|_{p/2}^{\delta}\leq c_{p,\delta,X}\sum_{k=1}%
^{n}\frac{1}{k^{1+2\delta/p}} \Biggl(\sum_{j=1}^{n}\sum_{i=j}^{\infty}%
c_{i} \Biggr)^{\delta} ,
\]
where $c_{p,\delta,X}$ is a positive constant depending depending only
on $p$,
$\delta$, $c$, $\kappa$ and $\Vert\eta_{0}\Vert_{p}$. Applying Theorem
\ref{directprop} for the martingale case, we then get the following corollary:
\begin{corollary}
\label{corARCH} Let $X=(X_{i})_{i\in{\mathbf{Z}}}$ be defined by
(\ref{defARCH}) and $S_{n}=\sum_{i=1}^{n}X_{i}$. Let $p>2$ and assume that
(\ref{condARCH}) is satisfied.\vspace*{8pt}

(1) If we assume that $\sum_{j\geq
n}c_{j}=O(n^{-b})$ for $b>1-2/p$, then, for any integer~$n$,
\[
{\mathbf{E}} \Bigl({\max_{1\leq j\leq n}}|S_{j}|^{p} \Bigr)\ll(n{\mathbf{E}%
}(X_{0}^{2}) )^{p/2}+n\bigl({\mathbf{E}}(|X_{0}|^{p})+b_{p,X}\bigr) ,
\]
where $b_{p,X}$ is a positive constant depending on $p$ and on the underlying
sequence~$X$, but not depending on $n$.

(2) If we assume $\sum_{j\geq
n}c_{j}=O(n^{-b})$ for $b>0$, then%
\[
\limsup_{n\rightarrow\infty}n^{-1/2}\Bigl\Vert{\max_{1\leq k\leq n}}
|S_{k}|\Bigr\Vert_{p}\leq a_{p}( {\mathbf{E}}(X_{0}^{2}))
^{1/2}<\infty,
\]
where $a_{p}$ is a constant depending only on $p$.
\end{corollary}

\subsection{Application to functions of linear processes}
\label{sectionFLP}

Let $(a_{i})_{i\in{{\mathbf{Z}}}}$ be a sequence of real numbers in
$\ell^{2}$
and $(\varepsilon_{i})_{i\in\mathbf{Z}}$ be a sequence of i.i.d. random
variables in ${\mathbf{L}}^{2}$. Define%
%
%
\begin{equation}\label{def2suite}%
X_{k}=h\biggl(\sum_{i\in{\mathbf{Z}}}a_{i}\varepsilon_{k-i}\biggr)-{\mathbf{E}%
}\biggl(h\biggl(\sum_{i\in{\mathbf{Z}}}a_{i}\varepsilon_{k-i}\biggr)\biggr) .\vadjust{\goodbreak}
\end{equation}
Denote by $w_{h}(\cdot,M)$ the modulus of continuity of the function $h$ on the
interval $[-M,M]$, that is,
\[
w_{h}(t,M)=\sup\{|h(x)-h(y)|,|x-y|\leq t,|x|\leq M,|y|\leq M\} .
\]

We shall establish the following result:
\begin{corollary}
\label{Thlin} Let $X=(X_{k})_{k\in{\mathbf{Z}}}$ be defined by
(\ref{def2suite}). Assume that $h$ is $\gamma$-H\"{o}lder on any
compact set,
with $w_{h}(t,M)\leq Ct^{\gamma}M^{\alpha}$, for some $C>0$, $\gamma\in\
]0,1]$
and $\alpha\geq0$. Let $p>2$ and assume that ${\mathbf{E}}(|\varepsilon
_{0}|^{2\vee(\alpha+\gamma)p})<\infty$.\vspace*{8pt}

(1) If $p$ is an even integer,
and for $\lambda>p/2-2$ if $p>4$ and $\lambda=0$ if $p=4$,
%
%
\begin{equation}\label{gammaholder}%
\sum_{i\geq1}i^{1-2/p}(\log i)^{\lambda}\Biggl(\sum_{j\geq i}a_{j}%
^{2}\Biggr)^{\gamma/2}<\infty,
\end{equation}
then, for any integer $n$,
\[
{\mathbf{E}} \Bigl({\max_{1\leq j\leq n}}|S_{j}|^{p} \Bigr)\ll({\mathbf{E}%
}(S_{n}^{2}) )^{p/2}+n\bigl({\mathbf{E}}(|X_{0}|^{p})+b_{p,X}\bigr) ,
\]
where $b_{p,X}$ is a positive constant depending on $p$ and on the underlying
sequence~$X$, but not depending on $n$.

(2) If for some $\eta>1$
%
%
\begin{equation} \label{cond2}%
\biggl(\sum_{j\geq i}a_{j}^{2} \biggr)^{\gamma/2}\ll i^{-\eta} ,
\end{equation}
then
\[
\limsup_{n\rightarrow\infty} n^{-1/2}\Bigl\Vert{\max_{1\leq k\leq n}}
|S_{k}|\Bigr\Vert_{p}\leq a_{p} \biggl(\sum_{k\in{\mathbf{Z}}}{\mathbf{E}%
}(X_{0}X_{k}) \biggr)^{1/2}<\infty,
\]
where $a_{p}$ is a constant depending only on $p$.
\end{corollary}
\begin{remark}
\label{rmklinear} The proof of the first item of the above result is
based on
Theorem~\ref{general}. Our proof reveals that an application of Theorem
\ref{directprop} would involve a more restrictive condition on $\lambda$,
namely $\lambda>p-3$.
\end{remark}

As a preliminary step in the proof of Corollary~\ref{Thlin} we state
the following proposition, which is a direct consequence of the proof
of Proposition 4.2 and of Theorem 4.2 in \citet{DedMerRio09},
page 988.
\begin{proposition}
\label{modulo} Let $(X_{i})_{i\in{\mathbf{Z}}}$ be as in Corollary~\ref{Thlin}
and $p\geq2$. Let $(\varepsilon_{i}^{\prime})_{i\in{\mathbf{Z}}}$ be an
independent copy of $(\varepsilon_{i})_{i\in{\mathbf{Z}}}$, and denote
$V_{0}=\sum_{i\geq0}a_{i}\varepsilon_{-i}$ and
\[
M_{1,i}=|V_{0}|\vee\biggl|\sum_{0\leq j<i}a_{j}\varepsilon_{-j}+\sum_{j\geq
i\geq0}a_{j}\varepsilon_{-j}^{\prime}\biggr|
\]
and
\[
\tilde{w}%
_{h}(i)= \biggl\|w_{h}\biggl(\biggl|\sum_{k\geq i}a_{k}\varepsilon_{-k}%
\biggr|,M_{1,i}\biggr) \biggr\|_{p}.
\]
Then for any $i,j\geq0$,
\[
\|{\mathbf{E}}_{0}(X_{i})\|_{p}\leq2\tilde{w}_{h}(i)
\]
and%
\begin{eqnarray*}
&&
\Vert{\mathbf{E}}_{0}(X_{i}X_{j+i})-{\mathbf{E}}(X_{i}X_{j+i})\Vert_{p/2}\\
&&\qquad\leq2\|X_{0}\|_{p}\min\bigl(\tilde{w}_{h}(i)+\tilde{w}_{h}(i+j),\tilde{w}%
_{h}(j-[j/2]) \bigr).
\end{eqnarray*}
Moreover, if $w_{h}(t,M)$ and $\varepsilon_{0}$ are as in Corollary
\ref{Thlin}, then for all $i\geq0$
\[
\tilde{w}_{h}(i)\leq K\biggl(\sum_{j\geq i}a_{i}^{2}\biggr)^{\gamma/2} ,
\]
where $K$ is a constant depending on $p$, $C$, $\alpha$, $\gamma$,
$\Vert\varepsilon_{0}\Vert_{2\vee(\alpha+\gamma)p}$ and on $\sum_{j\geq
0}%
a_{j}^{2}$.
\end{proposition}

\begin{pf*}{Proof of Corollary~\ref{Thlin}}
Notice first that by
Proposition~\ref{modulo} and relation~(\ref{cond2}), we get
%
%
\begin{equation} \label{Gordin}%
\sum_{k\geq1}\|{\mathbf{E}}_{0}(X_{k})\|_{p}<\infty.
\end{equation}
For any even integer $p\geq4$, applying Theorem~\ref{general} we
obtain, via relation (\ref{Gordin}) and the triangular inequality, that
\begin{eqnarray*}
&&{\mathbf{E}} \Bigl({\max_{1\leq j\leq n}}|S_{j}|^{p} \Bigr) \\
&&\qquad \ll
n \bigl(1+{\mathbf{E}}(|X_{0}|^{p})\bigr)+n \Biggl(\sum_{k=1}^{n}\frac{1}{k^{1+4/p(p-2)}
}({\mathbf{E}}(S_{k}^{2}))^{2/(p-2)} \Biggr)^{p(p-2)/4}\\
&&\qquad\quad{} +n \Biggl(\sum_{k=1}^{n}\frac{1}{k^{1+4/p(p-2)}}\|{\mathbf{E}}_{0}
(S_{k}%
^{2})-{\mathbf{E}}(S_{k}^{2})\|_{p/2}^{2/(p-2)} \Biggr)^{p(p-2)/4} .
\end{eqnarray*}
Notice now that relation (\ref{Gordin}) implies the so-called coboundary
decomposition; see Theorem 5.4 in \citet{HalHey80}; namely, there
is a
constant $K_{p,X}$ such that for all $n\geq1$
\[
S_{n}=M_{n}+R_{n}\qquad\mbox{with }\|R_{n}\|_{p}\leq2\sum_{k\geq1}\|{\mathbf
{E}%
}_{0}(X_{k})\|_{p}=K_{p,X} ,
\]
where $M_{n}$ is a martingale in ${\mathbf{L}}^{p}$ with stationary
differences. It easily follows that for all $1\leq k\leq n$,
\[
\frac{{\mathbf{E}}(S_{k}^{2})}{k}\leq4\frac{{\mathbf{E}}(S_{n}^{2})}{n}%
+2\frac{{\mathbf{E}}(R_{k}^{2})}{k}+4\frac{{\mathbf{E}}(R_{n}^{2})}{n}%
\leq4\frac{{\mathbf{E}}(S_{n}^{2})}{n}+6\frac{K_{p,X}}{k} .
\]
So,
\[
n \Biggl(\sum_{k=1}^{n}\frac{1}{k^{1+4/p(p-2)}}({\mathbf{E}}(S_{k}%
^{2}))^{2/(p-2)} \Biggr)^{p(p-2)/4}\ll({\mathbf{E}}(S_{n}^{2}))^{p/2}%
+K_{p,X}n.
\]
Therefore, the first part of the corollary follows if we prove that
(\ref{gammaholder}) implies that
%
%
\begin{equation}\label{cons1gammaholder*}%
\sum_{k\geq1}\frac{1}{k^{1+4/p(p-2)}}\|{\mathbf{E}}_{0}(S_{k}^{2}%
)-{\mathbf{E}}(S_{k}^{2})\|_{p/2}^{2/(p-2)}<\infty.
\end{equation}
Notice now that by H\"{o}lder's inequality (when $p > 4$) and (\ref
{varsquare}%
), in order to prove~(\ref{cons1gammaholder*}), it suffices to prove
that for
$\lambda>p/2-2$ if $p>4$ and $\lambda=0$ if $p=4$,
\[
\sum_{k\geq1}\frac{(\log k)^{\lambda}}{k^{1+2/p}}\sum_{i=1}^{k}\sum
_{j=0}^{k-i}\Vert{\mathbf{E}}_{0}(X_{i}X_{j+i})-{\mathbf{E}}(X_{i}%
X_{j+i})\Vert_{p/2}<\infty.
\]
The first part of Corollary~\ref{Thlin} is then a consequence of Proposition
\ref{modulo}. The second part follows in the same way by using Theorem
\ref{directprop}, simple computations and the fact that under (\ref{Gordin}),
$\lim_{n\rightarrow\infty}n^{-1}{\mathbf{E}}(S_{n})^{2}=\sum_{k\in
{\mathbf{Z}%
}}{\mathbf{E}}(X_{0}X_{k})$.~%
\end{pf*}

\subsection{Application to a stationary reversible Markov chain}
\label{sectionRMC}

First we want to mention that all our results can be formulated in the Markov
chain setting. We assume that $(\zeta_{n})_{n\in\mathbb{Z}}$ denotes a
stationary Markov chain defined on a probability space $(\Omega,\mathcal
{A}%
,{\mathbf{P}})$ with values in a measurable space $(E,\mathcal{E})$. The
marginal distribution and the transition kernel are denoted by $\pi
(A)={\mathbf{P}}(\zeta_{0}\in A)$ and $Q(\zeta_{0},A)={\mathbf{P}}(\zeta
_{1}\in A| \zeta_{0})$. In addition $Q$ denotes the operator acting via
$(Qf)(\zeta)=\int_{E}f(s)Q(\zeta,ds)$. Next, let $f$ be a function on
$E$ such
that $\int_{E}|f|^{p}\,d\pi<\infty$ and $\int_{E}f\,d\pi=0$. Denote by
$\mathcal{F}_{k}$ the $\sigma$-field generated by $\zeta_{i}$ with
$i\leq
k$, $X_{i}=f(\zeta_{i})$, and $S_{n}(f)=\sum_{i=1}^{n}X_{i}$. Notice
that any stationary sequence $(Y_{k})_{k\in\mathbb{Z}}$ can be viewed
as a
function of a Markov process $\zeta_{k}=(Y_{i};i\leq k)$, for the function
$g(\zeta_{k})=Y_{k}$.

The Markov chain is called reversible if $Q=Q^{\ast}$, where $Q^{\ast}$
is the
adjoint operator of $Q$. In this setting, an application of Theorem
\ref{stateven} gives the following estimate:
\begin{corollary}
\label{Rev}Let $(\zeta_{n})$ be a reversible Markov chain. For any even
integer $p\geq4$ and any positive integer $n$,
\begin{eqnarray*}
{\mathbf{E}} \Bigl({\max_{1\leq k\leq n}}|S_{k}(f)|^{p} \Bigr)&\ll& n\mathbf{E}%
(|f(\zeta_{1})|^{p})+n \Biggl(\sum_{k=1}^{n}\frac{1}{k^{1+1/p}}\Vert{\mathbf
{E}%
}_{0}(S_{k}(f))\Vert_{p} \Biggr)^{p}\\
&&{}+n \Biggl(\sum_{k=1}^{n}\frac{1}{k^{1+2/p}}\Vert{\mathbf{E}}_{0}
(S_{k}%
^{2}(f))-{\mathbf{E}}(S_{k}^{2}(f))\Vert_{p/2} \Biggr)^{p/2}\\
&&{} +n \Biggl(\sum
_{k=1}^{n}\frac{1}{k^{1+2/p}}{\mathbf{E}}(S_{k}^{2}(f)) \Biggr)^{p/2} .
\end{eqnarray*}
\end{corollary}

Moreover, using Theorem~\ref{pinteger}, we obtain:
\begin{corollary}
\label{Revgen} Let $(\zeta_{n})$ be a reversible Markov chain. For any real
number $p>4$ and any positive integer $n$,
\begin{eqnarray*}
{\mathbf{E}} \Bigl({\max_{1\leq k\leq n}}|S_{k}(f)|^{p} \Bigr)&\ll& n\mathbf{E}%
(|f(\zeta_{1})|^{p})\\
&&{ }+n \Biggl(\sum_{k=1}^{n}\frac{1}{k^{1+1/p}}\Vert{\mathbf{E}}_{0}(S_{k}%
^{2}(f))-{\mathbf{E}}(S_{k}^{2}(f))\Vert_{p/2}^{1/2} \Biggr)^{p}\\
&&{} +n \Biggl(\sum
_{k=1}^{n}\frac{1}{k^{1+1/p}}{\|}S_{k}(f)\|_{2} \Biggr)^{p} .
\end{eqnarray*}
\end{corollary}

This corollary is also valid for any real $2 < p \leq4$. For this range,
however, according to the Comment~\ref{commentgreaterdelta}, Theorem
\ref{directprop} (resp., Corollary~\ref{Rev}) gives a better bound for
$p\in\ ]2,4[$ (resp., for $p=4$).

For a particular example, let $E=[-1,1]$, and let $\upsilon$ be a symmetric
atomless law on $E$. The transition probabilities are defined by
\[
Q(x,A)=(1-|x|)\delta_{x}(A)+|x|\upsilon(A) ,
\]
where $\delta_{x}$ denotes the Dirac measure. Assume that $\theta=\int
_{E}|x|^{-1}\upsilon(dx)<\infty$. Then there is a unique invariant measure
\[
\pi(dx)=\theta^{-1}|x|^{-1}\upsilon(dx),
\]
and the stationary Markov chain $(\zeta_{i})_{i}$ is reversible and
positively recurrent.

Assume the following assumption on the measure $\upsilon$: there exists a
positive constant $c$ such that for any $x\in[0,1]$,
%
%
\begin{equation}\label{hyp2upsilon}%
\frac{d\upsilon}{dx}(x)\leq cx^{p/2-1}\bigl(\log(1+1/x)\bigr)^{-\lambda}
\qquad\mbox{for some
$\lambda>0$}.
\end{equation}
As an application of Corollary~\ref{Revgen} we shall establish:
\begin{corollary}
\label{cor2MarkovChain} Let $p>2$ be a real number, and let
$f(-x)=-f(x)$ for
any $x\in E$. Assume that $|f(x)|\leq C|x|^{1/2}$ for any $x$ in $E$
and a
positive constant~$C$.

(1) Assume in addition that (\ref{hyp2upsilon})
is satisfied for $\lambda>p$. Then for any integer $n$,
%
%
\begin{eqnarray} \label{revineq}%
{\mathbf{E}} \Bigl({\max_{1\leq k\leq n}}|S_{k}(f)|^{p} \Bigr)&\ll& n^{p/2}%
\biggl(\int_{0}^{1}f^{2}(x)x^{-2}\upsilon(dx) \biggr)^{p/2}\nonumber\\[-8pt]\\[-8pt]
&&{}+n\bigl({\mathbf{E}%
}(|f(\zeta_{0})|^{p})+b_{p,\lambda,c,C}\bigr) ,\nonumber
\end{eqnarray}
where $b_{p,\lambda,c,C}$ is a positive constant depending on $p$,
$\lambda$,
$c$ and $C$.

(2) Assume now (\ref{hyp2upsilon}) is relaxed to
\[
\frac{d\upsilon}{dx}(x)\leq cx^{a}\qquad\mbox{for some $a>0$.}%
\]
Then
\[
\limsup_{n\rightarrow\infty} n^{-1/2}\Bigl\Vert{\max_{1\leq k\leq n}}
|S_{k}(f)|\Bigr\Vert_{p}\leq a_{p} \biggl(\int_{0}^{1}f^{2}(x)x^{-2}%
\upsilon(dx) \biggr)^{1/2} ,
\]
where $a_{p}$ is a constant depending only on $p$.
\end{corollary}

Notice that this example of reversible Markov chain has been considered
by Rio
[(\citeyear{Rio09}), Section 4] under a slightly more stringent
condition on the measure
than (\ref{hyp2upsilon}). Corollary~\ref{cor2MarkovChain} then extends
Proposition 4.1(b) in \citet{Rio09} to all real numbers $p>2$.
\begin{pf*}{Proof of Corollary~\ref{cor2MarkovChain}}
To get this result we shall apply Corollary~\ref{Revgen}. We start by noticing
that $f$ being an odd function, we have
%
%
\begin{equation}\label{operator}%
{\mathbf{E}}(f(\zeta_{k})|\zeta_{0})=(1-|\zeta_{0}|)^{k}f(\zeta
_{0})\qquad\mbox{a.s.}
\end{equation}
Therefore, for any $j\geq0$,
\[
{\mathbf{E}}(X_{0}X_{j})={\mathbf{E}}(f(\zeta_{0}){\mathbf{E}}(f(\zeta
_{j})|\zeta_{0}))=\theta^{-1}\int_{E}f^{2}(x)(1-|x|)^{j}|x|^{-1}%
\upsilon(dx) .
\]
Then, by inequality (\ref{var}),
%
%
\begin{eqnarray}\label{estimatevar}
{\mathbf{E}}(S_{k}^{2}(f)) & \leq & 2k\theta^{-1} \biggl(\int_{0}^{1}%
f^{2}(x)x^{-1}\upsilon(dx)+2\sum_{j=1}^{k-1}\int
_{0}^{1}f^{2}(x)(1-x)^{j}%
x^{-1}\upsilon(dx) \biggr)\hspace*{-14pt}\nonumber\\[-8pt]\\[-8pt]
& \leq & 2k\theta^{-1} \biggl(\int_{0}^{1}f^{2}(x)x^{-1}\upsilon(dx)+2\int_{0}%
^{1}f^{2}(x)x^{-2}\upsilon(dx) \biggr) .\nonumber
\end{eqnarray}
Next, we give an upper bound of the quantity $\Vert{\mathbf
{E}}_{0}(S_{n}%
^{2}(f))-{\mathbf{E}}(S_{n}^{2}(f))\Vert_{p/2}$. By using the fact that for
any positive $k$, $\pi Q^{k}=\pi$ we have
\begin{eqnarray*}
&&
\Vert{\mathbf{E}}_{0}(S_{n}^{2}(f))-{\mathbf{E}}(S_{n}^{2}(f))\Vert
_{p/2}\\
&&\qquad\leq\sum_{k=1}^{n} \Biggl(\int_{E} \Biggl|(\delta_{x}Q^{k}-\pi) \Biggl(f^{2}%
+2f\sum_{k=1}^{n-k}Q^{k}f \Biggr) \Biggr|^{p/2}\pi(dx) \Biggr)^{2/p} ;
\end{eqnarray*}
see also inequality (4.12) in \citet{Rio09}. Now, by using the relation
(\ref{operator}), one can prove that for any $x\in E$,
\[
\Biggl(f^{2}+2f\sum_{k=1}^{n}Q^{k}f \Biggr)(x)=f^{2}(x)
\bigl(1+2\bigl(1-(1-|x|)^{n}%
\bigr)(|x|^{-1}-1) \bigr) ;
\]
see the computations in \citet{Rio09} leading to his relation
(4.13). Then,
since $|f(x)|\leq C|x|^{1/2}$, it follows that
\[
\sup_{x\in E} \Biggl|f^{2}(x)+2f(x)\sum_{k=1}^{n}%
Q^{k}f(x) \Biggr|\leq2C^{2} .
\]
Therefore, for any $p>2$,
%
%
\begin{equation}\label{maje0sn2}%
\Vert{\mathbf{E}}_{0}(S_{n}^{2}(f))-{\mathbf{E}}(S_{n}^{2}(f))\Vert
_{p/2}%
\leq4C^{2}\sum_{k=1}^{n} \biggl(\int_{E}\Vert Q^{k}(x,\cdot)-\pi(\cdot)\Vert
\pi(dx) \biggr)^{2/p} ,\hspace*{-28pt}
\end{equation}
where $\Vert\mu(\cdot)\Vert$ denotes the total variation of the signed measure
$\mu$; see also inequality (4.15) in \citet{Rio09}. We estimate
next the
coefficients of absolute regularity $\beta_{n}$ as defined in
(\ref{boundbetacor2}). Let $a=p/2-1$. We shall prove now that under
(\ref{hyp2upsilon}), there exists a positive constant $K$ depending on $a$,
$c$ and $\lambda$, such that
%
%
\begin{equation} \label{boundbetacor2}%
2\beta_{n}:=\int_{E}\Vert Q^{n}(x,\cdot)-\pi(\cdot)\Vert\pi(dx)\leq
Kn^{-a}(\log n)^{-\lambda} .
\end{equation}
Notice first that by [Lemma 2, page 75, in \citet{DouMasRio94}], we
have that
%
%
\begin{equation}\label{relbetatau}%
\beta_{n}\leq3\int_{E}(1-|x|)^{[n/2]}\pi(dx) .
\end{equation}
Let $k\geq2$ be an integer. Clearly, for any $\alpha\in\ ]0,1[$,
%
%
\begin{eqnarray}\label{calculmomentT}%
\int_{0}^{1}(1-x)^{k}\pi(dx)&\leq& c\int_{0}^{k^{-\alpha}}(1-x)^{k}x^{a-1}
\bigl(\log(1+1/x)\bigr)^{-\lambda}\,dx\nonumber\\[-8pt]\\[-8pt]
&&{}+c\int_{k^{-\alpha}}^{1}(1-x)^{k}x^{a-1}%
\bigl(\log(1+1/x)\bigr)^{-\lambda}\,dx .\nonumber
\end{eqnarray}
Notice now that
\[
\int_{0}^{k^{-\alpha}}(1-x)^{k}x^{a-1}\bigl(\log(1+1/x)\bigr)^{-\lambda}\,dx\leq
(\alpha\log k)^{-\lambda}\int_{0}^{1}(1-x)^{k}x^{a-1}\,dx .
\]
Hence, by the properties of the Beta and Gamma functions,
%
%
\begin{equation}\label{p1boundbeta}\quad
\lim_{k\rightarrow\infty}k^{a}(\log k)^{\lambda}\int_{0}^{k^{-\alpha}%
}(1-x)^{k}x^{a-1}\bigl(\log(1+1/x)\bigr)^{-\lambda}\,dx\leq\alpha^{-\lambda}a\Gamma
(a) .
\end{equation}
On the other hand, we have that
\[
\int_{k^{-\alpha}}^{1}(1-x)^{k}x^{a-1}\bigl(\log(1+1/x)\bigr)^{-\lambda}\,dx\leq
(\log2)^{-\lambda}(1-k^{-\alpha})^{k}\int_{0}^{1}x^{a-1}\,dx ,
\]
and then, since $\alpha<1$, we easily obtain
%
%
\begin{equation}\label{p2boundbeta}%
\lim_{k\rightarrow\infty}k^{a}(\log k)^{\lambda}\int_{k^{-\alpha}}%
^{1}(1-x)^{k}x^{a-1}\bigl(\log(1+1/x)\bigr)^{-\lambda}\,dx=0 .
\end{equation}
Starting from (\ref{relbetatau}) and taking into account (\ref
{calculmomentT}%
), (\ref{p1boundbeta}) and (\ref{p2boundbeta}), (\ref{boundbetacor2}) follows.
Then, by using inequality (\ref{maje0sn2}) combined with
(\ref{boundbetacor2}), we derive that
\[
\Vert{\mathbf{E}}_{0}(S_{n}^{2}(f))-{\mathbf{E}}(S_{n}^{2}(f))\Vert
_{p/2}\leq
C^{2}K_{p,c,\lambda}n^{2/p}(\log n)^{-2\lambda/p} ,
\]
where $K_{p,c,\lambda}$ is a positive constant depending on $p,c$ and
$\lambda$. Therefore, for any $\delta\in\ ]0,1]$, there exists a positive
constant $b_{p,\lambda,c,C}$ depending on $p$, $\lambda$, $c$ and $C$, such
that
%
%
\begin{equation}\label{convergencelog}%
\sum_{k=1}^{n}\frac{1}{k^{1+1/p}}\Vert{\mathbf{E}}_{0}(S_{k}^{2}%
(f))-{\mathbf{E}}(S_{k}^{2}(f))\Vert_{p/2}^{1/2}\leq b_{p,\lambda
,c,C}\qquad\mbox{for $\lambda>p$} .
\end{equation}
Considering estimates (\ref{estimatevar}) and (\ref{convergencelog}), the
first part of Corollary~\ref{cor2MarkovChain} follows from an
application of
Corollary~\ref{Revgen} (taking also into account the comment after its
statement for $2<p\leq4$). To get the second part, we also apply Corollary
\ref{Revgen} but with the upper bound $\beta_{n}=O(n^{-a})$ instead of
(\ref{boundbetacor2}).
\end{pf*}

\section{Application to density estimation}
\label{densitysection}

In this section, we estimate the ${\mathbf{L}}^{p}$-integrated risk for
$p\geq4$, for the kernel estimator of the unknown marginal density $f$
of a
stationary sequence $(Y_{i})_{i\geq0}$.

Applying our Theorem~\ref{directprop}, we shall show that if the coefficients
of dependence $(\beta_{2,Y}(k))_{k\geq1}$ (see Definition~\ref{beta}
below) of the sequence $(Y_{i})_{i\in{\mathbf{Z}}}$ satisfy $\beta
_{2,Y}(k)=O(n^{-a})$ for $a>p-1$, then the bound of the ${\mathbf
{L}}^{p}%
$-norm of the random term of the risk is of the same order of magnitude
as the
optimal one obtained in \citet{BreHub79} in the independence
setting (see their Corollary 2), provided that the density is bounded
and the
kernel $K$ satisfies Assumption~\ref{assumpAp} below.

\renewcommand{\theAssumption}{A$_{p}$}
\begin{Assumption}\label{assumpAp}
$K$ is a BV
(bounded variation) function such that
\[
\int_{\mathbf{R}}|K(u)|\,du<\infty
\quad\mbox{and}\quad\int_{\mathbf{R}}|K(u)|^{p}%
\,du<\infty.
\]
\end{Assumption}
\begin{definition}
\label{beta} Let $(Y_{i})_{i\in{\mathbf{Z}}}$ be a stationary sequence
of real-valued random variables, and let ${\mathcal{F}}_{0}=\sigma
(Y_{i},i\leq0)$. For
any positive $i$ and $j$, define the random variables
\[
b({\mathcal{F}}_{0},i,j)=\sup_{(s,t)\in{\mathbf{R}}^{2}}\bigl|\mathbf{P}%
(Y_{i} \leq t,Y_{j} \leq s|{\mathcal{F}}_{0})-\mathbf{P}(Y_{i} \leq
t,Y_{j} \leq s)\bigr|.
\]
Define now the coefficient
\[
\beta_{2,Y}(k)=\sup_{i\geq j\geq k}{\mathbf{E}}(b({\mathcal
{F}}_{0},i,j)) .
\]
\end{definition}
\begin{proposition}
\label{applidensity} Let $p\geq4$ and $K$ be any real function satisfying
Assumption~\ref{assumpAp}. Let $(Y_{i})_{i\geq0}$ be a stationary sequence
with unknown marginal density $f$ such that $\Vert f\Vert_{\infty
}<\infty$.
Define
\[
X_{k,n}(x)=K\bigl(h_{n}^{-1}(x-Y_{k})\bigr)
\quad\mbox{and}\quad f_{n}(x)=\frac{1}{nh_{n}}%
\sum_{k=1}^{n}X_{k,n}(x) ,
\]
where $(h_{n})_{n\geq1}$ is a sequence of positive real
numbers.\vspace*{8pt}

(1)
Assume that there exists a positive constant $c$ such that for some
$\eta>0$
and all $n\geq1$,
%
%
\begin{equation}\label{condbetathm}%
\beta_{2,Y}(n)\leq cn^{-(p-1+\eta)} .
\end{equation}
Then there exist positive constants $C_{1}$ and $C_{2}$ depending on $p$,
$\eta$ and $c$ such that for any positive integer $n$,
%
%
\begin{eqnarray}\label{resdensity}
&&
{\mathbf{E}}\int_{\mathbf{R}}|f_{n}(x)-{\mathbf{E}}(f_{n}(x))|^{p}\,dx\nonumber\\
&&\qquad\leq
C_{1}(nh_{n})^{-p/2}\|dK\|^{p/2}\|f\|_{\infty}^{p/2-1} \biggl(\int_{\mathbf{R}
}|K(u)|\,du \biggr)^{p/2}\nonumber\\[-8pt]\\[-8pt]
&&\qquad\quad{}+C_{2}(nh_{n})^{1-p} \biggl(\int_{\mathbf{R}}|K(u)|^{p}\,du+\|dK\| \int
_{\mathbf{R}}|K(u)|^{p-1}\,du\nonumber\\
&&\qquad\quad\hspace*{134.3pt}{}+\|dK\|^{2}\int_{\mathbf{R}}|K(u)|^{p-2}%
\,du \biggr) ,\nonumber
\end{eqnarray}
where $\|dK\|$ is the total variation norm of the measure $dK$.

(2)
Assume now, in addition, that $nh_{n}\rightarrow\infty$. Then we also
have the
following asymptotic result:%
\begin{eqnarray*}
&&\limsup_{n\rightarrow\infty}(nh_{n})^{-p/2}{\mathbf{E}}\int
_{\mathbf{R}}|f_{n}(x)-{\mathbf{E}}(f_{n}(x))|^{p}\,dx\\
&&\qquad\leq a_{p,c,\eta
}\|dK\|^{p/2}\|f\|_{\infty}^{p/2-1} \biggl(\int_{\mathbf{R}}|K(u)|\,du \biggr)^{p/2},
\end{eqnarray*}
where $a_{p,c,\eta}$ is a constant depending only on $p$, $c$ and
$\eta$.
\end{proposition}

The bound obtained in Proposition~\ref{applidensity} can be also
compared to
the one obtained in Theorem 3.3 in \citet{Vie97} under the
assumption that
the strong $\beta$-mixing coefficients in the sense of \citet
{VolRoz59} of the sequence $(Y_{i})_{i\in{\mathbf{Z}}}$, denoted by
$\beta
_{\infty}(k)$, satisfy $\sum_{k\geq1}k^{p-2}\beta_{\infty}(k)<\infty$. Our
condition is then comparable to the one imposed by \citet{Vie97}
but less
restrictive in the sense that many processes are such that the sequence
$\beta_{2,Y}(n)$ tends to zero as $n\rightarrow\infty$ which is not the case
for $\beta_{\infty}(n)$; see the examples given in \citet
{DedPri07}. Notice also that, for $p=2$, inequality (\ref{resdensity}) is
proven in Dedecker and Prieur [(\citeyear{DedPri05}), Proposition~3] under the summability of
the weak $\beta$-dependence coefficients. In addition, when $p=3$,
under the
condition $\sum_{k\geq1}k\beta_{2,Y}(k)<\infty$, it can be proven that
inequality (\ref{resdensity}) still holds by applying the Rosenthal-type
inequality in \citet{Ded01} [see also Proposition 3 in \citet
{DedPri07} for a version of the
Dedecker's inequality in case of stationary
sequences]; however, his inequality does not lead to (\ref{resdensity}) when
$p>3$.

When the random variables $Y_{i}$ are a function of an i.i.d. sequence,
namely, $Y_{i}=g(\ldots,\varepsilon_{i-1},\varepsilon_{i})$ where $g$
is a
measurable function, and $(\varepsilon_{n},n\in{\mathbf{Z}})$ are
i.i.d. random
variables, Wu, Huang and Huang [(\citeyear{WuHuaHua10}), Theorem 1]
obtained an upper bound of
similar order as in (\ref{resdensity}) (for $p>1$) under conditions
imposed to
the so-called ``predictive dependent measures.'' When, in addition, the variable
$\varepsilon_{0}$ has a density with bounded derivatives up to order
$2$, the
method is especially effective for short memory linear processes with
independent innovations; see Section 4.1 in \citet{WuHuaHua10}. Our
Proposition~\ref{applidensity} complements the above cited results
since the
coefficients $\beta_{2,Y}(n)$ can be\vspace*{1pt} estimated without assuming that
$\varepsilon_{0}$ has a density. For instance,\vspace*{1pt} we obtain upper bound
(\ref{resdensity}) for $Y_{i}=\sum_{k\geq0}2^{-k-1}\varepsilon_{i-k}$,
and the
$\varepsilon_{k}$'s are i.i.d. Bernoulli random variables with
parameter $1/2$; see Section~6.1 of \citet{DedPri07} for
computations of
$\beta_{2,Y}(n)$. In addition, our Proposition~\ref{applidensity} applies
even for situations when the variables $Y_{i}$ are not assumed to be a
function of an i.i.d. sequence. We refer, for instance, to
Dedecker and Prieur
[(\citeyear{DedPri09}), Theorem 3.1] who gave an upper bound of the
coefficients $\beta
_{2,Y}(n)$ of the Markov chain associated to an intermittent map.

If we assume now that $f$ has a derivative of order $s$, where $s\geq1$
is an
integer and that the following bound holds for the bias term:
%
%
\begin{equation}\label{biais}%
\int_{\mathbf{R}}|f(x)-{\mathbf{E}}(f_{n}(x))|^{p}\,dx\leq
Mh_{n}^{sp}\bigl\Vert
f^{(s)}\bigr\Vert^{p}_{p} ,
\end{equation}
where $M$ is a constant depending on the kernel $K$, then the choice of
$(nh_{n})^{p/2}\times h_{n}^{sp}=O(1)$ leads to the following estimate:
%
%
\begin{equation}\label{erreurtotale}%
{\mathbf{E}}\int_{\mathbf{R}}|f_{n}(x)-f(x)|^{p}\,dx=O\bigl(n^{-sp/(2s+1)}\bigr) .
\end{equation}
We mention that (\ref{biais}) holds for any Parzen kernel of order $s$; see
Section 4 in \citet{BreHub79}. We also mention that if we only
assume that $\sum_{k\geq1}k^{p-2}\beta_{2,Y}(k)<\infty$ instead of
(\ref{condbetathm}) in Proposition~\ref{applidensity}, then inequality
(\ref{resdensity}) is valid with $(nh_{n})^{1-p}n^{\varepsilon}$ (for any
$\varepsilon>0$) replacing $(nh_{n})^{1-p}$ in the second term of the
right-hand side. In this situation, bound (\ref{biais}) combined with a
choice of $h_{n}$ of order $n^{-1/(1+2s)}$ still leads to estimate
(\ref{erreurtotale}).
\begin{pf*}{Proof of Proposition~\ref{applidensity}}
Setting $X_{i,n}(x)=K((x-Y_{i})/h_{n})-\break{\mathbf
{E}}(K((x-Y_{i})/h_{n}))$, we
have that
%
%
\begin{equation}\label{trivialbounddensity}%
{\mathbf{E}}\int_{\mathbf{R}}|f_{n}(x)-{\mathbf{E}}(f_{n}(x))|^{p}%
\,dx\leq(nh_{n})^{-p}\int_{\mathbf{R}}{\mathbf{E}} \Biggl|\sum_{i=1}^{n}%
X_{i,n}(x) \Biggr|^{p}\,dx .
\end{equation}
Starting from (\ref{trivialbounddensity}) and applying Proposition
\ref{consdirect2} to the stationary sequence $(X_{i,n}(x))_{i\in\mathbf{Z}}$,
Proposition~\ref{applidensity} follows provided we establish the following
bounds (in what follows $C$ is a positive constant which may vary from
line to
line, and that may depend on $p$, $c$ and $\eta$ but not on $n$):
%
%
\begin{equation}\label{estdensityb0}%
\int_{\mathbf{R}}{\mathbf{E}}|X_{1,n}(x)|^{p}\,dx\leq2^{p+1}h_{n}\int
_{\mathbf{R}}|K(u)|^{p}\,du ,
\end{equation}
%
%
\begin{eqnarray} \label{estdensityb1}
\nonumber\\[-14pt]
&&\int_{\mathbf{R}}\Biggl(\sum_{j=0}^{n-1}|{\mathbf{E}}(X_{0,n}(x)X_{j,n}%
(x))| \Biggr)^{p/2}\,dx\nonumber\\[-8pt]\\[-8pt]
&&\qquad\leq Ch_{n}^{p/2}\|dK\|^{p/2}\|f\|_{\infty}^{(p/2)-1}%
\biggl(\int_{\mathbf{R}}|K(u)|\,du \biggr)^{p/2} ,\nonumber
\end{eqnarray}
and that for $\varepsilon>0$ small enough,
%
%
\begin{eqnarray}\label{estdensityb3}%
&&\sum_{j=1}^{n}j^{p-2+\varepsilon}\int_{\mathbf{R}}\Vert
X_{0,n}(x){\mathbf{E}%
}_{0}(X_{j,n}(x))\Vert_{p/2}^{p/2}\,dx\nonumber\\[-8pt]\\[-8pt]
&&\qquad\leq Ch_{n}\|dK\| \int_{\mathbf{R}%
}|K(u)|^{p-1}\,du \nonumber
\end{eqnarray}
and
%
%
\begin{eqnarray}\label{estdensityb4}\qquad
&&\sum_{j=1}^{n}j^{p-2+\varepsilon}\sup_{i\geq j}  \int_{\mathbf{R}}%
\Vert{\mathbf{E}}_{0}(X_{i,n}(x)X_{j,n}(x))-{\mathbf{E}}(X_{i,n}%
(x)X_{j,n}(x))\Vert_{p/2}^{p/2}\,dx\nonumber\\[-8pt]\\[-8pt]
&&\qquad \leq Ch_{n}\|dK\|^{2} \int_{\mathbf{R}}|K(u)|^{p-2}\,du .\nonumber
\end{eqnarray}
In what follows, we shall prove these bounds. Notice first that
\[
\int_{\mathbf{R}}{\mathbf{E}}|X_{1,n}(x)|^{p}\,dx\leq2^{p+1}\int_{\mathbf
{R}%
}\int_{\mathbf{R}}\bigl|K\bigl((x-y)h_{n}^{-1}\bigr)\bigr|^{p}f(y)\,dx\,dy
\]
proving (\ref{estdensityb0}) by the change of variables $u=(x-y)h_{n}^{-1}$.
To prove (\ref{estdensityb1}), we first apply item (1) of Lemma~\ref{covbeta}
implying that%
\[
\sum_{j=0}^{n-1}|{\mathbf{E}}(X_{0,n}(x)X_{j,n}(x))|\leq\|dK\| {\mathbf
{E}%
} \bigl(\widetilde{b}({\mathcal{F}}_{0},n)
\bigl|K\bigl((x-Y_{0})/h_{n}\bigr)\bigr| \bigr) ,
\]
where $\widetilde{b}({\mathcal{F}}_{0},n)=\sum_{j=0}^{n-1}b({\mathcal
{F}}%
_{0},j,j)$. An application of H\"{o}lder's inequality as done in
Viennet [(\citeyear{Vie97}), page 474], then gives
\begin{eqnarray*}
&&\int_{\mathbf{R}}\Biggl(\sum_{j=0}^{n-1}|{\mathbf{E}}(X_{0,n}(x)X_{j,n}%
(x))| \Biggr)^{p/2}\,dx\\
&&\qquad\leq h_{n}^{p/2}\Vert f\Vert_{\infty}^{p/2-1}{\mathbf{E}%
} (\widetilde{b}({\mathcal{F}}_{0},n))^{p/2} \biggl(\int_{\mathbf{R}%
}|K(u)|\,du \biggr)^{p/2} .
\end{eqnarray*}
This proves (\ref{estdensityb1}) since ${\mathbf{E}} (\widetilde
{b}({\mathcal{F}}_{0},n))^{p/2}\leq C\sum_{k=1}^{n}k^{p-2}\beta
_{2,Y}(k)$ and
\[
\sum_{k=1}^{n}k^{p-2}\beta_{2,Y}(k)=O(1)
\]
by condition (\ref{condbetathm}).

We turn now to the proof of (\ref{estdensityb3}). With this aim we
notice that
\[
\Vert X_{0,n}(x){\mathbf{E}}_{0}(X_{j,n}(x))\Vert_{p/2}^{p/2}={\mathbf
{E}%
}(Z_{0}(x){\mathbf{E}}_{0}(X_{j,n}(x)))={\mathbf
{E}}(Z_{0}(x)X_{j,n}(x)) ,
\]
where $Z_{0}(x)=|X_{0,n}|^{p/2}|{\mathbf{E}}_{0}(X_{j,n}(x))|^{p/2-1}%
\operatorname{sign} ({\mathbf{E}}_{0}(X_{j,n}(x)) )$. Then, by using
item (1)
of Lemma~\ref{covbeta}, we derive that
%
%
\begin{eqnarray}\label{used1lemmacovbeta}%
\Vert X_{0,n}(x){\mathbf{E}}_{0}(X_{j,n}(x))\Vert_{p/2}^{p/2}&=&
\operatorname{Cov}%
\bigl(Z_{0}(x),K\bigl((x-Y_{j})/h_{n}\bigr) \bigr)\nonumber\\[-8pt]\\[-8pt]
&\leq&\|dK\| {\mathbf{E}}%
(b({\mathcal{F}}_{0},j,j) |Z_{0}(x)| ) .\nonumber
\end{eqnarray}
Notice now that by using the elementary inequality: $x^{\alpha
}y^{1-\alpha
}\leq x+y$ valid for $\alpha\in[0,1]$ and nonnegative $x$ and $y$, we
get that $|Z_{0}(x)|\leq(|X_{0,n}(x)|+|{\mathbf{E}}_{0}(X_{j,n}%
(x))| )^{p-1}$. Therefore, some computations involving Jensen's inequality
lead to
%
%
\begin{equation}\label{ine1Z0}%
\int_{\mathbf{R}}|Z_{0}(x)|\,dx\leq4^{p}h_{n}\int_{\mathbf
{R}}|K(u)|^{p-1}\,du ,
\end{equation}
where $\|dK\|$ is the total variation norm of the measure $dK$.
Starting from
(\ref{used1lemmacovbeta}), we end the proof of (\ref{estdensityb3}) by taking
into account (\ref{ine1Z0}) and the fact that
\[
\sum_{j=1}^{n}j^{p-2+\varepsilon}{\mathbf{E}}(b({\mathcal{F}}_{0}%
,j,j))\leq\sum_{j=1}^{n}j^{p-2+\varepsilon}\beta_{2,Y}(j)
\]
is convergent by condition (\ref{condbetathm}) for any $\varepsilon<\eta$.\vadjust{\goodbreak}

It remains to prove (\ref{estdensityb4}). We first write that
\[
\Vert{\mathbf{E}}_{0}(X_{i,n}(x)X_{j,n}(x))-{\mathbf{E}}(X_{i,n}%
(x)X_{j,n}(x))\Vert_{p/2}^{p/2}={\mathbf{E}}\bigl((Z_{0}^{\prime})^{(0)}%
(x)X_{i,n}(x)X_{j,n}(x)\bigr) ,
\]
where the notation $X^{(0)}$ stands for $X^{(0)}=X-{\mathbf{E}}(X)$ and
\[
Z_{0}^{\prime}(x)=|{\mathbf{E}}_{0}(B_{i,j}(x))|^{p/2-1}\operatorname
{sign}%
({\mathbf{E}}_{0}(B_{i,j}(x)) )
\]
with $B_{i,j}(x)=X_{i,n}(x)X_{j,n}(x)-{\mathbf{E}}(X_{i,n}(x)X_{j,n}(x))$.
Since the variables $X_{i,n}(x)$ and $X_{j,n}(x)$ are centered, an application
of item (2) of Lemma~\ref{covbeta} then gives
\begin{eqnarray*}
&& \Vert{\mathbf{E}}_{0}(X_{i,n}(x)X_{j,n}(x))-{\mathbf{E}}(X_{i,n}%
(x)X_{j,n}(x))\Vert_{p/2}^{p/2}\\
&&\qquad \leq\|dK\|^{2} {\mathbf{E}} \bigl(|Z_{0}(x)|\bigl(b({\mathcal{F}}_{0}%
,i,i)+b({\mathcal{F}}_{0},j,j)+b({\mathcal{F}}_{0},i,j)\bigr) \bigr).
\end{eqnarray*}
Notice now that since $p/2-1\geq1$, we can easily get
\[
\int_{\mathbf{R}}|Z_{0}^{\prime}(x)|\,dx\leq c_{p}h_{n}\int_{\mathbf{R}%
}|K(u)|^{p-2}\,du ,
\]
where $c_{p}$ is a positive constant depending on $p$. In addition,
\[
\sum_{j=1}^{n}j^{p-2+\varepsilon}\sup_{i\geq j}{\mathbf{E}}\bigl(b({\mathcal
{F}%
}_{0},i,i)+b({\mathcal{F}}_{0},j,j)+b({\mathcal{F}}_{0},i,j)\bigr)\leq3\sum
_{j=1}^{n}j^{p-2+\varepsilon}\beta_{2,Y}(j) ,
\]
which is convergent by condition (\ref{condbetathm}) for any
$\varepsilon
<\eta$. Then (\ref{estdensityb4}) holds and so does the proposition.
\end{pf*}

\begin{appendix}\label{app}
\section*{Appendix}

This section is devoted to some technical lemmas. The next lemma gives estimates
for terms of the type ${\mathbf{E}}(X_{0}^{u}X_{1}^{p-u})$.
\begin{lemma}
\label{basic} Let $p$ and $u$ be real numbers such that $0\leq u\leq
p-2$. Let
$X_{0}$ and $X_{1}$ be two positive identically distributed random variables.
With the notation $a^{p}={\mathbf{E}}(X_{0}^{p})$, ${\mathbf{E}}_{0}%
(X_{1})={\mathbf{E}}(X_{1}|X_{0})$ the following estimates hold:
%
%
\begin{equation} \label{b1basic}%
{\mathbf{E}}(X_{0}^{u}X_{1}^{p-u})\leq a^{p-2u/(p-2)}\|{\mathbf{E}}_{0}%
(X_{1}^{2})\|_{p/2}^{u/(p-2)}
\end{equation}
and
%
%
\begin{equation}\label{b2basic}%
{\mathbf{E}}(X_{0}^{p-1}X_{1})\leq a^{p-1}\|{\mathbf{E}}_{0}(X_{1}%
^{2})\|_{p/2}^{1/2} .
\end{equation}
\end{lemma}
\begin{pf}
Inequality (\ref{b1basic})
is trivial for $u=0$. To prove it for $u=p-2$, it suffices to write that
${\mathbf{E}}(X_{0}^{p-2}X_{1}^{2})={\mathbf{E}}(X_{0}^{p-2}{\mathbf{E}}
_{0}(X_{1}^{2}))$, and then to use H\"{o}lder's inequality.

We prove now inequality (\ref{b1basic}) for $0<u<p-2$. Select
$x=(p/2-1)/u=(p-2)/2u$. Notice that\vadjust{\goodbreak} $2x>1$ and $p-u-1/x>0$ since $u<p-2$.
Then, since the variables are identically distributed,
\begin{eqnarray*}
{\mathbf{E}}(X_{0}^{u}X_{1}^{p-u})&=&{\mathbf{E}}(X_{0}^{u}X_{1}^{1/x}%
X_{1}^{p-u-1/x})\leq\|X_{0}^{u}X_{1}^{1/x}\|_{2x}\|X_{1}^{p-u-1/x}%
\|_{2x/(2x-1)}\\
&\leq&({\mathbf{E}}(X_{0}^{p-2}X_{1}^{2}) )^{u/(p-2)}(a^{p}%
)^{1-u/(p-2)} .
\end{eqnarray*}
Now, again by H\"{o}lder's inequality applied with $x=p/(p-2)$ and
$1-1/x=2/p$,
\begin{eqnarray*}
{\mathbf{E}}(X_{0}^{p-2}X_{1}^{2})&=&{\mathbf{E}}(X_{0}^{p-2}{\mathbf{E}}%
_{0}(X_{1}^{2}))\leq({\mathbf{E}}(X_{0}^{p}) )^{(p-2)/p}%
({\mathbf{E}}({\mathbf{E}}_{0}(X_{1}^{2}))^{p/2} )^{2/p}\\
&=&a^{(p-2)}\|{\mathbf{E}}_{0}(X_{1}^{2})\|_{p/2} .
\end{eqnarray*}
Overall
\begin{eqnarray*}
{\mathbf{E}}(X_{0}^{u}X_{1}^{p-u}) & \leq & a^{u}\|{\mathbf{E}}_{0}(X_{1}%
^{2})\|_{p/2}^{u/(p-2)}(a^{p})^{1-u/(p-2)}\\
& = &a^{p-2u/(p-2)}\|{\mathbf{E}}_{0}(X_{1}^{2})\|_{p/2}^{u/(p-2)}
\end{eqnarray*}
ending the proof of inequality (\ref{b1basic}).

To prove inequality (\ref{b2basic}), we use H\"{o}lder's inequality which
entails that
\[
{\mathbf{E}}(X_{0}^{p-1}X_{1})\leq{\mathbf{E}}(X_{0}^{p-1}{\mathbf{E}}%
_{0}^{1/2}(X_{1}^{2}))\leq a^{p-1}\|{\mathbf{E}}_{0}(X_{1}^{2})\|_{p/2}%
^{1/2} .
\]
\upqed\end{pf}

The next lemma gives covariance-type inequalities in terms of beta
coefficients as
defined in Definition~\ref{beta}.
\begin{lemma}
\label{covbeta}Let $Z$ be a ${\mathcal{F}}_{0}$-measurable real-valued random
variable, and let $h$ and $g$ be two BV functions. Denote by $\|dh\|$
(resp.,
$\|dg\|$) the total variation norm of the measure $dh$ (resp., $dg$). Denote
$Z^{(0)}=Z-{\mathbf{E}}(Z)$, $h^{(0)}(Y_{i})=h(Y_{i})-{\mathbf{E}}(h(Y_{i}))$
and $g^{(0)}(Y_{j})=g(Y_{j})-{\mathbf{E}}(g(Y_{j}))$. Define the random
variables $b({\mathcal{F}}_{0},i,j)$ as in Definition~\ref{beta}. Then:

\begin{longlist}[(2)]
\item[(1)] $ |{\mathbf{E}} (Z^{(0)}h^{(0)}(Y_{i}%
) ) |= |{\operatorname{Cov}}(Z,h(Y_{i})) |\leq\|dh\| {\mathbf{E}%
} (|Z|b({\mathcal{F}}_{0},i,i) )$.

\item[(2)] $ |{\mathbf{E}} (Z^{(0)}h^{(0)}(Y_{i})g^{(0)}(Y_{j}%
) ) |\leq\|dh\|\|dg\| {\mathbf{E}} (|Z|(b({\mathcal{F}}%
_{0},i,i)+b({\mathcal{F}}_{0},j,j)+\break b({\mathcal{F}}_{0}, i,j) )$.
\end{longlist}
\end{lemma}
\begin{pf}
Item (1) has been proven by
\citet{DedPri05} [see item (2) of their Proposition 1]. Item (2) needs a
proof. We first notice that
\[
h^{(0)}(X)g^{(0)}(Y)=\iint\bigl(\mathbf{1}_{X\leq t}-F_{X}%
(t)\bigr)\bigl(\mathbf{1}_{Y\leq s}-F_{Y}(s)\bigr)\,dh(t)\,dg(s) .
\]
Therefore
\begin{eqnarray*}
&& {\mathbf{E}} \bigl(Z^{(0)}h^{(0)}(Y_{i})g^{(0)}(Y_{j}) \bigr)\\
&&\qquad ={\mathbf{E}} \biggl(Z\iint\bigl(\mathbf{1}_{Y_{i}\leq
t}^{(0)}\mathbf{1}_{Y_{j}\leq s}^{(0)}-{\mathbf{E}} \bigl(\mathbf
{1}_{Y_{i}\leq
t}^{(0)}\mathbf{1}_{Y_{j}\leq s}^{(0)} \bigr) \bigr)\,dh(t)\,dg(s) \biggr)\\
&&\qquad ={\mathbf{E}} \biggl(Z\iint{\mathbf{E}} \bigl(\mathbf{1}%
_{Y_{i}\leq t}^{(0)}\mathbf{1}_{Y_{j}\leq s}^{(0)}-{\mathbf{E}}%
\bigl(\mathbf{1}_{Y_{i}\leq t}^{(0)}\mathbf{1}_{Y_{j}\leq s}^{(0)}%
\bigr)|{\mathcal{F}}_{0} \bigr)\,dh(t)\,dg(s) \biggr) ,
\end{eqnarray*}
which proves item (2) by noticing that
\begin{eqnarray*}
\bigl|{\mathbf{E}} \bigl(\mathbf{1}_{Y_{i}\leq t}^{(0)}\mathbf{1}_{Y_{j}\leq
s}^{(0)}-{\mathbf{E}} \bigl(\mathbf{1}_{Y_{i}\leq t}^{(0)}\mathbf
{1}_{Y_{j}\leq
s}^{(0)} \bigr)|{\mathcal{F}}_{0} \bigr)\bigr|
\leq b({\mathcal{F}}_{0}%
,i,i)+b({\mathcal{F}}_{0},j,j)+b({\mathcal{F}}_{0},i,j) .\quad
\end{eqnarray*}
\upqed\end{pf}

The next lemma gives inequalities for $|x+y|^{p}$ for different ranges
of $p$,
where $p\geq2$ is a real number.
\begin{lemma}
\label{cross}
(1) Let $x$ and $y$ be two real numbers and $2\leq p\leq3$. Then%
%
%
\begin{equation}\label{Rio}%
|x+y|^{p}\leq|x|^{p}+|y|^{p}+p|x|^{p-1}\operatorname{sign}(x)y+\frac
{p(p-1)}%
{2}|x|^{p-2}y^{2} .
\end{equation}

(2) Let $x$ and $y$ be two real numbers and $3<p\leq4$. Then%
%
%
\begin{eqnarray}\label{Rio34}%
|x+y|^{p}&\leq&|x|^{p}+|y|^{p}+p|x|^{p-1}\operatorname{sign}(x)y+\frac
{p(p-1)}%
{2}|x|^{p-2}y^{2}\nonumber\\[-8pt]\\[-8pt]
&&{}+\frac{2p}{(p-2)}|x||y|^{p-1} .\nonumber
\end{eqnarray}

(3) Let $x$ and $y$ be two positive real numbers and $p\geq1$
any real
number. Then%
%
%
\begin{equation}\label{int}%
(x+y)^{p}\leq x^{p}+y^{p}+4^{p}(x^{p-1}y+xy^{p-1}) .
\end{equation}

(4) Let $x$ and $y$ be two real numbers and $p$ an even positive integer.
Then%
%
%
\begin{equation}\label{evenint}\quad
(x+y)^{p}\leq
x^{p}+y^{p}+p(x^{p-1}y+xy^{p-1})+2^{p}(x^{2}y^{p-2}+x^{p-2}%
y^{2}) .
\end{equation}
\end{lemma}
\begin{pf}
Inequality (\ref{Rio}) was
established in Rio [(\citeyear{Rio09}), Relation (3.3)] by using Taylor's
expansion with
integral rest for evaluating the difference $|x+y|^{p}-|x|^{p}$. To prove
inequality (\ref{Rio34}), we also use the Taylor integral formula of
order $2$
that gives
%
%
\begin{eqnarray}\label{order2}%
|x+y|^{p}-|x|^{p}&=&p|x|^{p-1}\operatorname
{sign}(x)y+C_{p}^{2}|x|^{p-2}y^{2}\nonumber\\[-8pt]\\[-8pt]
&&{}+2C_{p}^{2}y^{2}\int_{0}^{1}(1-t)(|x+ty|^{p-2}-|x|^{p-2})\,dt
,\nonumber
\end{eqnarray}
where $C_{p}^{2}=p(p-1)/2$. Notice now that, for $3<p\leq4$,
\[
|x+ty|^{p-2}\leq\frac{x^{2}+2|x||ty|+y^{2}}{(|x|+|ty|)^{4-p}}\leq
|x|^{p-2}+2|x||ty|^{p-3}+|ty|^{p-2} .
\]
Hence
\begin{eqnarray*}
&& 2C_{p}^{2}y^{2}\int_{0}^{1}(1-t)(|x+ty|^{p-2}-|x|^{p-2})\,dt\\
&&\qquad \leq2C_{p}^{2}|y|^{p}\int_{0}^{1}(1-t)t^{p-2}\,dt+4C_{p}^{2}|x||y|^{p-1}
\int_{0}^{1}(1-t)t^{p-3}\,dt\\
&&\qquad =2|y|^{p}C_{p}^{2}\frac{\Gamma(p-1)}{\Gamma(p+1)}+4|x||y|^{p-1}C_{p}%
^{2}\frac{\Gamma(p-2)}{\Gamma(p)}\\
&&\qquad=|y|^{p}+\frac{2p}{(p-2)}|x||y|^{p-1} .
\end{eqnarray*}
Starting from (\ref{order2}), 
inequality (\ref{Rio34}) follows.

Inequality (\ref{int}) was observed by \citet{Sha95}, page 957. We shall
establish now (\ref{evenint}). We start by noticing that for any $a$,
$b$ two
positive real numbers and $2\leq k\leq p-2$ we have
%
%
\begin{equation}\label{ab}%
a^{p-k}b^{k}\leq\max\bigl(a^{p}(b/a)^{p-2},a^{p}(b/a)^{2}\bigr)\leq a^{2}b^{p-2}%
+a^{p-2}b^{2} .
\end{equation}

Now for $p$ an even positive integer and $x$ and $y$ two real
numbers, the Newton binomial formula gives
\begin{eqnarray*}
(x+y)^{p} & = & x^{p}+y^{p}+p(x^{p-1}y+xy^{p-1})+\sum_{k=2}^{p-2}C_{p}%
^{k}x^{p-k}y^{k}\\
& \leq & x^{p}+y^{p}+p(x^{p-1}y+xy^{p-1})+\sum_{k=2}^{p-2}C_{p}^{k}%
|x|^{p-k}|y|^{k} .
\end{eqnarray*}
Whence, by (\ref{ab}) and the fact that $\sum_{k=0}^{p}C_{p}^{k}=2^{p}$,
inequality (\ref{evenint}) follows.~%
\end{pf}
\begin{lemma}
\label{lmasubadd} Let $(V_{i})_{i \geq0}$ be a sequence of nonnegative
numbers such that $V_{0} = 0$ and for all $i,j \geq0$,
%
%
\begin{equation}
\label{condsubadd}
V_{i+j} \leq C ( V_{i} + V_{j}) ,
\end{equation}
where $C \geq1$ is a constant not depending on $i$ and $j$. Then:

\begin{longlist}[(3)]
\item[(1)] For any integer $r \geq1$, any integer $n$ satisfying
$2^{r-1} \leq n <
2^{r}$ and any real $q \geq0$,
\[
\sum_{i=0}^{r-1} \frac{1}{2^{iq}} V_{2^{i}} \leq C 2^{q+2} ( 2^{q+1} -1)^{-1}
\sum_{k= 1}^{n} \frac{V_{k}}{k^{1+q}} .
\]

\item[(2)] For any positive integers $k$ and $m$ and any real $q>0$,
\[
\sum_{j=1}^{k}\frac{1}{j^{q}}V_{jm}\leq2^{q+1}Cq^{-1}m^{q-1}\sum_{\ell
=1}%
^{m}\frac{1}{(\ell+m)^{q}}V_{\ell}+2Cq^{-1}m^{q-1}\sum_{\ell
=m+1}^{km}\frac
{1}{\ell^{q}}V_{\ell} .
\]

\item[(3)] Let $0<\delta\leq\gamma\leq1$. Then for any real $q\geq0$,
\[
\Biggl( \sum_{k=1}^{n}\frac{1}{k^{1+q\gamma}}V_{k}^{\gamma}\Biggr)
^{1/\gamma}\leq2^{1/\delta-1/\gamma}C^{(\gamma-\delta)/\delta}\Biggl(
\sum_{k=1}^{n}\frac{1}{k^{1+q\delta}}V_{k}^{\delta}\Biggr) ^{1/\delta} .
\]
\end{longlist}
\end{lemma}
\begin{remark}
If $(V_{i})_{i \geq0}$ satisfies (\ref{condsubadd}) with $C=1$, then the
sequence is said to be subadditive.\vadjust{\goodbreak}
\end{remark}
\begin{pf*}{Proof of Lemma~\ref{lmasubadd}}
Condition (\ref{condsubadd}) implies that for any integer $k$ and any
integer $0\leq j\leq k$,
%
%
\begin{equation}\label{condsubadd2}%
V_{k}\leq C(V_{j}+V_{k-j})\quad\mbox{and then that}\quad(k+1)V_{k}\leq2C\sum
_{j=1}%
^{k}V_{j} .
\end{equation}
Therefore for $2^{r-1}\leq n<2^{r}$,
\[
\sum_{i=0}^{r-1}\frac{1}{2^{iq}}V_{2^{i}}\leq2C\sum_{j=1}^{2^{r-1}}V_{j}
\sum_{i\dvtx2^{i}\geq j}\frac{1}{2^{i(q+1)}}
\]
proving item (1). To prove item (2), using again (\ref{condsubadd2}), it suffices
to notice that
\begin{eqnarray*}
\sum_{j=1}^{k}\frac{1}{j^{q}}V_{jm} & \leq & 2C\sum_{j=1}^{k}\frac
{1}{(1+jm)j^{q}}\sum_{\ell=1}^{jm}V_{\ell}\\
& \leq & 2Cm^{-1}(2m)^{q}\sum_{j=1}^{k}j^{-q-1}\sum_{\ell=1}^{m}\frac{1}%
{(\ell+m)^{q}}V_{\ell}\\
&&{}+2Cq^{-1}m^{q-1}\sum_{\ell=m+1}^{km}\frac{1}{\ell
^{q}%
}V_{\ell} .
\end{eqnarray*}
To prove item (3), we first notice that (\ref{condsubadd}) entails that
\[
V_{i+j}^{\gamma}\leq C^{\gamma}(V_{i}^{\gamma}+V_{j}^{\gamma})
\quad\mbox{and then
that}\quad(k+1)V_{k}^{\gamma}\leq2C^{\gamma}\sum_{j=1}^{k}V_{j}^{\gamma} .
\]
Then for any real $q\geq0$,
%
%
\begin{equation}\label{ineitem3}%
k^{-q(\gamma-\delta)}V_{k}^{\gamma-\delta}\leq2^{1-\delta/\gamma}%
C^{\gamma-\delta} \Biggl(\sum_{j=1}^{k}j^{-(1+q\gamma)}V_{j}^{\gamma
} \Biggr)^{1-\delta/\gamma} .
\end{equation}
Writing that $k^{-(1+q\gamma)}V_{k}^{\gamma}=(k^{-(1+q\delta
)}V_{k}^{\delta
})(k^{-q(\gamma-\delta)}V_{k}^{\gamma-\delta})$ and using (\ref{ineitem3}),
the following inequality holds:
\[
\sum_{k=1}^{n}\frac{1}{k^{1+q\gamma}}V_{k}^{\gamma}\leq2^{1-\delta
/\gamma
}C^{\gamma-\delta} \Biggl(\sum_{k=1}^{n}k^{-(1+q\delta)}V_{k}^{\delta
}\Biggr) \Biggl(\sum_{j=1}^{n}j^{-(1+q\gamma)}V_{j}^{\gamma}\Biggr)^{1-\delta
/\gamma}
\]
proving item (3).
\end{pf*}
\end{appendix}

\section*{Acknowledgments}

The authors are indebted to the referees for carefully reading the
manuscript and for helpful comments that improved the presentation of
the paper.


%

\printaddresses

\end{document}